\documentclass{article}

\usepackage{amsfonts}
\usepackage{amssymb}
\usepackage{epsfig}
\usepackage{amsthm}

\newcommand{\N}{{\mathbb N}}
\newcommand{\Z}{{\mathbb Z}}
\newcommand{\ga}{$\Gamma(G,{\cal A})$ }
\renewcommand{\O}{{\cal O}}
\newcommand{\h}{{\hat h}}

\newtheorem{lemma}{Lemma}[section]
\newtheorem{cor}{Corollary}
\newtheorem{thm}{Theorem}
\newtheorem{prop}{Proposition}

\newcounter{remarks}
\newcommand{\remark}{\par \refstepcounter{remarks}{\noindent \bf Remark \arabic{remarks}. }}

\begin{document}
\begin{center}
{\large SOME PROPERTIES OF SUBSETS OF HYPERBOLIC GROUPS

\vspace{.5cm}
 Ashot Minasyan}

\vspace{.3cm}
Department of Mathematics

Vanderbilt University

Nashville, TN 37240, USA

aminasyan@gmail.com
\end{center}

\begin{abstract}
We present some results about  quasiconvex subgroups
of infinite index and their products. After that we
extend the standard notion of a subgroup commensurator to an
arbitrary subset of a group, and generalize some of the previously
known results.
\end{abstract}

\section{\large Introduction}

Assume that $G$ is a $\delta$-hyperbolic group for some $\delta \ge 0$ and \ga is its Cayley
graph corresponding to some finite symmetrized (i.e. ${\cal A}={\cal A}^{-1}$) generating set $\cal A$.
\ga is a proper geodesic metric space; a subset $Q \subseteq G$ is said to be
$\varepsilon$-{\it quasiconvex}, if any geodesic connecting two elements from $Q$
belongs to a closed $\varepsilon$-neighborhood ${\cal O}_{\varepsilon}(Q)$ of $Q$ in $\Gamma(G,{\cal A})$
for some $\varepsilon \ge 0$. $Q$ will be called {\it quasiconvex} if there exists $\varepsilon > 0$ for which
it is $\varepsilon$-quasiconvex.

In \cite{Gromov} Gromov shows that the notion of quasiconvexity in a hyperbolic group
does not depend on the choice of a finite generating set .

If $A,B \subset G$, $x \in G$ we define $A^B = \{bab^{-1}~|~a \in A,b \in B\}$, $A^x = xAx^{-1}$.
 For any $x \in G$ the set $x^G=\{x\}^G$ is the conjugacy class of the element $x$.

If $x \in G$, $o(x)$ will denote the
order of the element $x$ in $G$, $|x|_G$ -- the length of a shortest representation of $x$ in
terms of the generators from $\cal A$. $\langle x \rangle$ will be the cyclic subgroup of $G$
generated by $x$ (sometimes, if $o(x)=n$, we will write $\langle x \rangle_n$).
 $1_G$ will denote the identity element of $G$. For a subgroup
$H \le G$, $|G:H|$ will be the index of $H$ in $G$.

\begin{thm} \label{prop1} Let $H_1,H_2,\dots,H_s$ be quasiconvex subgroups of a hyperbolic group $G$.
Let $K$ be an arbitrary subgroup of $G$. Then the following two
conditions are equivalent:

(a) $|K:(K\cap H_j^g)| = \infty$ for every $j \in \{1,2,\dots,s\}$ and every $g\in G$;

(b) there exists an element of infinite order $x \in K$
such that the intersection $\langle x \rangle_\infty \cap (H_1^G
\cup H_2^G \cup \dots \cup H_s^G)$ is trivial.
\end{thm}

Let $G_1,G_2,\dots,G_n$ be quasiconvex subgroups of $G$,
$f_0,f_1,\dots,$ $f_n \in G$, \\ $n \in \N \cup \{ 0 \}$. Using the same terminology as in \cite{hyp},
the set
$$P=f_0G_1f_1G_2 \cdot \dots\cdot f_{n-1}G_nf_n = \{ f_0g_1f_1 \cdot \dots \cdot g_nf_n \in G~|~g_i
\in G_i,~i=1,\dots,n\} $$ will be called a {\it quasiconvex product}.
\\ The quasiconvex subgroups $G_i$, $i=1,2,\dots,n$, are {\it members} of the product $P$.

Let $U=\bigcup_{k=1}^q P_k$ be a finite union of quasiconvex products $P_k$, $k=1,\dots,q$.
A subgroup $H \le G$ will be called a { \it member} of $U$, by definition, if $H$ is a member of $P_k$
for some $1 \le k \le q$. For any such set $U$ we fix its representation as a finite union of quasiconvex
products and fix its members.

\begin{thm} \label{prop2} Assume that $U$ is a finite union of quasiconvex products in
a hyperbolic group $G$ and the subgroups $H_1,H_2,\dots,H_s$ are all the members of $U$.
If $K$ is a subgroup of $G$ and $K \subseteq U$ then
for some $g \in G$ and $j\in \{1,2,\dots,s\}$ one has $|K:(K\cap H_j^g)|<\infty$.
\end{thm}

We will say that a finite union of quasiconvex products has {\it infinite
index} in $G$ if each of its members has infinite index in $G$. As an immediate consequence of Theorem\ref{prop2} applied to the case when $K=G$ we achieve

\begin{cor} Let $G$ be a hyperbolic group and $U$ be a finite union of quasiconvex products of
infinite index in $G$. Then $U$ is a proper subset of $G$, i.e. $G \neq U$.

\end{cor}
%
%

In section \ref{concepts} we generalize the definition of a subgroup commensurator. More precisely,
to any subset $A \subset G$ we correspond a subgroup $Comm_G(A) \le G$. In sections \ref{known} and
\ref{generalization} we list several known results about quasiconvex subgroups and their commensurators
and then extend them to more general settings.

A group is said to be non-elementary if it is not virtually cyclic.
In section \ref{conjsection} we investigate some properties of infinite conjugacy classes and prove

\begin{thm} \label{conjugacythm} Let $G$ be a non-elementary hyperbolic group and $A$ be
 a finite union of conjugacy classes in $G$. If the subset $A$ is infinite then it is not quasiconvex.
\end{thm}

\vspace{.3cm}
{\noindent \large \bf Acknowledgements}

The author would like to thank his advisor Professor A. Yu. Ol'shanskii
for his support and helpful comments.

\section{\large Some concepts and definitions} {\label {concepts}}

Suppose $G$ is an arbitrary group and $2^G$ is the set of all its subsets.
Below we establish some auxiliary relations on $2^G$. Assume $A,B \subseteq G$.

\vspace{.15cm}
{\bf \underline{Definition.}}  We will write
$B \preceq A$ if there exist elements $x_1,\dots,x_n \in G$ such that
$$B \subset Ax_1 \cup Ax_2 \cup \dots \cup Ax_n~.$$
If $G \preceq A$, the subset $A$ will be called {\it quasidense}.

\vspace{.15cm}
Obviously, the relation "$\preceq$" is transitive and reflexive.

Any subgroup $H$ of finite index in $G$ is quasidense; the complement \\
$H^{(c)}=G \backslash H$ in this case is also quasidense  (if $H\neq G$) since it contains a left coset modulo $H$,
and a shift (left or right) of a quasidense subset is quasidense.

On the other hand, if $H \le G$ and $|G:H| = \infty$, the set of elements of $H$ is not quasidense in
$G$. There is $y \in G \backslash H$, hence for any $x \in G$ either $x \in G \backslash H$ or
$xy \in G \backslash H$, thus $G = H^{(c)} \cup H^{(c)}y^{-1}$, i.e. $H^{(c)}$ is still quasidense.

\vspace{.15cm}
{\bf \underline{Definition.}} $A$ and $B$ will be called {\it equivalent} if $A \preceq B$ and $B \preceq A$.
In this case we will use the notation $A \approx B$.

\vspace{.15cm}
It is easy to check that "$\approx$" is an equivalence relation on $2^G$. Now, let $[A]$ denote
the equivalence class of a subset $A \subseteq G$ and let $\cal M$ be the set of all such equivalence
classes. Evidently, the relation "$\preceq$" induces a partial order on $\cal M$: $[A],[B] \in {\cal M}$,
$[A] \le [B]$ if and only if $A \preceq B$.

The group $G$ acts on $\cal M$ as follows: $g \in G$, $[A] \in {\cal M}$, then
$g \circ [A] = [gA]$. Indeed, the verification of the group action axioms is straightforward:

1. If $g,h \in G$, $[A] \in {\cal M}$ then $(gh) \circ [A] = g \circ (h \circ [A])$;

2. If $1_G \in G$ is the identity element and $[A] \in {\cal M}$ then $1_G \circ [A] = [A]$. \\
This action is well defined because if $A \approx B$, $g \in G$, then $gA \approx gB$.

If $A \subseteq G$, the stabilizer of $[A]$ under this action is the subgroup
 $$St_G([A]) = \{g \in G~|~g \circ [A] = [A]\}.$$

\vspace{.15cm}
{\bf \underline{Definition.}} For a given subset $A$ of the group $G$ the subgroup $St_G([A])$ will be called
{\it commensurator} of $A$ in $G$ and denoted $Comm_G(A)$. In other words,
$$Comm_G(A) = \{g \in G~|~gA \approx A\}~.$$

\vspace{.15cm}
Thus, to an arbitrary subset $A$ of the group we corresponded a subgroup in $G$.
Now, let's list some

\vspace{.15cm}
{\bf \noindent Properties of $Comm_G(A)$:}

1) If $card(A) < \infty$ or $A$ is quasidense then $Comm_G(A)=G$ (because any two finite non-empty subsets are equivalent
and a left shift of a  quasidense subset is quasidense);

2) If $A,B \subseteq G$ and $A \approx B$ then $Comm_G(A)=Comm_G(B)$;

3) The commensurator of $A \subset G$ contains (as its subgroups) the normalizer of $A$
~~$N_G(A) = \{g \in G~|~ gAg^{-1} = A\}$ and  the stabilizer under the action of the group $G$ on itself by left
multiplication $St_G(A) = \{g \in G~|~gA=A\}$.

4) For any $h \in G$ ~~$Comm_G(hA)=hComm_G(A)h^{-1}$.

\vspace{.15cm}
\begin{lemma} \label{sgpcompar}  Let $A,B$ be subgroups of $G$. Then
$A \preceq B$ if and only if the index $|A:(A \cap B)|$ is finite.
\end{lemma}

\underline{Proof.} The sufficiency is trivial. To prove the necessity,
suppose there exist $y_j \in G$, $j=1,2,\dots,m$, such that $A \subset By_1 \cup \dots \cup By_m$.
Without loss of generality we can assume that $A \cap By_j \neq \emptyset$ for every $j=1,2,\dots,m$.
Then for each $j=1,2,\dots,m$, there are $a_j \in A,b_j \in B$ such that $y_j=b_ja_j$. Hence
$By_j = Ba_j$ for all $j$, and therefore
$$A = \bigcup_{j=1}^m By_j \cap A = \bigcup_{j=1}^m (Ba_j \cap A) = \bigcup_{j=1}^m (B \cap A)a_j,$$
i.e. $|A:(B \cap A)| < \infty$. $\square$

\vspace{.15cm}
For a subgroup $H \le G$ the standard notion of the commensurator (virtual normalizer) subgroup of $H$ is given by
$$VN_G(H) = \{g \in G~|~|H:(H \cap gHg^{-1})| < \infty,|gHg^{-1}:(H \cap gHg^{-1})| < \infty\}~.$$
Now we are going to show that our new definition is just a generalization of it:

\vspace{.15cm}
\remark  \label{rem1} If $H$ is a subgroup of the group $G$ then $Comm_G(H) = VN_G(H)$.

\vspace{.15cm}
Indeed, let $g \in VN_G(H)$. Then, by definition, $$H \preceq (H \cap gHg^{-1}) \preceq gHg^{-1} \preceq gH~~\mbox{and}~~
gH \preceq gHg^{-1} \preceq (H \cap gHg^{-1}) \preceq H~,$$ thus $H \approx gH$ and $g \in Comm_G(H)$. So,
$VN_G(H) \subseteq Comm_G(H)$.

Now, suppose $g \in Comm_G(H)$, implying $H \approx gH$ but $gH \approx gHg^{-1}$, hence $H \preceq gHg^{-1}$
and $gHg^{-1} \preceq H$. By Lemma \ref{sgpcompar}, $g \in VN_G(H)$. Therefore $VN_G(H) = Comm_G(H)$.

If the group $G$ is finitely generated, one  can fix a  finite symmetrized generating set
$\cal A$ and define the word metric $d(\cdot,\cdot)$ corresponding to $\cal A$
in the standard way: first for every $g \in G$ we define
$l_{\cal A}(g)$ to be the length of a shortest word in $\cal A$ representing $g$; second, for any $x,y \in G$ we set
$d(x,y) = l_{\cal A}(x^{-1}y)$. Now, for arbitrary two subsets $A,B \subseteq G$ one can establish
$$h(A,B) = \inf\{\varepsilon >0~|~A \subset {\cal O}_\varepsilon(B), B \subset {\cal O}_\varepsilon(A)\}~\mbox{--}$$
the Hausdorff distance between $A$ and $B$ (${\cal O}_\varepsilon(B)$ is the closed $\varepsilon$-neighborhood
of $B$ in $G$). Where an infinum over the empty set is defined to be positive infinity.
\\ In this case, for any $A,B \subseteq G$ we observe that $B \preceq A$ if and only if there exists $c > 0$
such that $B \subset {\cal O}_c (A)$, and, therefore, $A \approx B$ if and only if $h(A,B) < \infty$.

Now we investigate the special case, when $G$ is $\delta$-hyperbolic (for the definition of a hyperbolic group
see section \ref{preliminaries}) and, therefore, finitely generated.
We will need the following

\vspace{.15cm}
\remark  \label{rem2} (\cite[Remark 4, lemma 2.1]{hyp}) Let $Q,A,B\subseteq G$ be quasiconvex subsets, $g\in G$.
Then (a)~the left shift $gQ=\{ gx~|~x\in Q\}$ is quasiconvex ;
(b)~the right shift $Qg=\{ xg~|~x\in Q\}$ is quasiconvex;
(c)$A \cup B$ is quasiconvex.

Therefore, a left coset modulo a quasiconvex subgroup and a conjugate subgroup to it
are quasiconvex in $G$.

\vspace{.15cm}
\remark  \label{rem2.5} Let the group $G$ be hyperbolic.

{\bf 1)} Suppose a subset $A \subset G$ is quasiconvex and $A \approx B$ for some $B \subset G$. Then $B$ is
also quasiconvex.

Indeed, as we saw above,
there exist $c_1,c_2 \ge 0$ such that $B \subset {\cal O}_{c_1}(A)$ and $A \subset {\cal O}_{c_2}(B)$.
Consider arbitrary $x,y \in B$ and a geodesic segment $[x,y]$ connecting them. Then
$$x,y \in {\cal O}_{c_1}(A) = \bigcup_{g\in G,|g|_G \le c_1} Ag~,$$ which is $\varepsilon$-quasiconvex
by Remark \ref{rem2} for some $\varepsilon \ge 0$. Therefore,

$$[x,y] \subset {\cal O}_{c_1+\varepsilon}(A) \subset {\cal O}_{c_1+\varepsilon+c_2}(B)$$ implying that
$B$ is $(c_1+c_2+\varepsilon)$-quasiconvex.

{\bf 2)} A subset $Q$ of the group $G$ (or of the Cayley
graph \ga) is quasidense if and only if there exists $\alpha \ge 0$ such that
for every $x \in G$ (or \ga) the distance $d(x,Q) = \inf\{d(x,y)~|~y \in Q \}$ is at most $ \alpha$,
i.e. $G \subseteq \O_{\alpha}(Q)$.

Indeed, if $G=Qg_1 \cup Qg_2 \cup \dots \cup Qg_n$, where $g_i \in G$, $i=1,2,\dots,n$,
Denote \\ $\alpha = max\{|g_i|_G~:~1\le i \le n\}$. Then for any $x \in G$,
there are $i \in \{1,\dots,n\}$ and $y \in Q$ with $x=yg_i$, hence
$d(y,x) = |g_i|_G \le \alpha$.

For demonstrating the sufficiency, let $\{g_1,g_2,\dots,g_n\}$ be the set of all elements in $G$ of
length at most $\alpha$. Then for every $x \in G$ there exists $y \in Q$ with $d(y,x) = |y^{-1}x|_G \le \alpha$;
hence, $y^{-1}x = g_i$ for some $i \in \{1,2,\dots,n\}$. Thus, $x = yg_i \in Qg_i$.

{\bf 3)} A quasidense subset $Q \subseteq G$ is quasiconvex.

This is an immediate consequence of the part 2).

\vspace{.15cm}

\section{\large Preliminaries} \label{preliminaries}

Assume $\bigl(X,d(\cdot,\cdot)\bigr)$ is a proper geodesic metric space.
If $Q \subset X$, $N \ge 0$, the closed $N$-neighborhood of $Q$ will be denoted by
$${\cal O}_N (Q) \stackrel{def}{=} \{x\in X~|~d(x,Q) \le N \}~.$$

If $x,y,w \in X$, then the number
$$(x|y)_w \stackrel{def}{=} \frac12 \Bigl(d(x,w)+d(y,w)-d(x,y) \Bigr)$$
is called the {\it Gromov product} of $x$ and $y$ with respect to $w$.

Let $abc$ be a geodesic triangle in the space $X$ and $[a,b]$, $[b,c]$, $[a,c]$ be its sides between the corresponding vertices.
There exist "special" points $O_a \in [b,c]$,\\ $O_b \in [a,c]$, $O_c \in [a,b]$
with the properties:
$d(a,O_b) = d(a,O_c) = \alpha$, $d(b,O_a) = $ $=d(b,O_c) = \beta$, $d(c,O_a) = d(c,O_b) = \gamma$. From a corresponding
system of linear equations one can find that $\alpha = (b|c)_a$, $\beta = (a|c)_b$, $\gamma = (a|b)_c$. Two points
$O \in [a,b]$ and $O' \in [a,c]$ are called $a$-{\it equidistant} if $d(a,O) = d(a,O') \le \alpha$.
The triangle $abc$ is said to be $\delta$-{\it thin} if for any two points $O,O'$ lying on its sides and
equidistant from one of its vertices, $d(O,O') \le \delta$ holds (Figure 0).

\begin{figure}
\begin{center}
\input{pic0.tex}
\end{center}
\begin{center}
Figure 0
\end{center}

\end{figure}

A geodesic $n$-gon in the space $X$ is said to be $\delta$-{\it slim} if each of its sides belongs to a
closed $\delta$-neighborhood of the union of the others.

We assume the following equivalent definitions of hyperbolicity of the space $X$
to be known to the reader (see \cite{Ghys},\cite{Mihalik}): \\
$1^{\circ}.$ There exists $\delta \ge 0$ such that for any four points $x,y,z,w \in X$ their Gromov products satisfy
$$(x|y)_w \ge min\{(x|z)_w,(y|z)_w\} - \delta~;$$
$2^{\circ}.$ All triangles in $X$ are {\cal $\delta$-thin} for some $\delta \ge 0$;\\
$3^{\circ}.$ All triangles in $X$ are $\delta$-slim for some $\delta \ge 0$.

Now, suppose $G$ is finitely generated group with a fixed finite symmetrized generating set $\cal A$.
One can define $d(\cdot,\cdot)$ to be the usual left-invariant metric on the Cayley graph of the group $G$
corresponding to $\cal A$. Then the Cayley graph \ga becomes a proper geodesic metric space.
$G$ is called hyperbolic if \ga is a hyperbolic metric space. It is easy to show that this definition does not depend on
the choice of the finite generating set $\cal A$ in $G$, thus hyperbolicity is a group-theoretical property.
It is well known that free groups of finite rank are hyperbolic.

Further on we will assume that \ga meets $1^{\circ},2^{\circ}$ and $3^{\circ}$ for a fixed (sufficiently large)
$\delta \ge 0$.

For any two points $x,y \in \Gamma(G,{\cal A})$ we fix a geodesic path between them and denote it by $[x,y]$.



Let $p$ be a path in the Cayley graph of $G$. Further on by $p_-$, $p_{+}$ we will denote the startpoint and
the endpoint of $p$, by $||p||$ -- its length; $lab(p)$, as usual, will mean the word in the alphabet $\cal A$
written on $p$. $elem(p) \in G$ will denote the element of the group $G$ represented by the word $lab(p)$.

A path $q$ is called $(\lambda,c)$-{ \it quasigeodesic} if there exist $0<\lambda \le 1$, $c \ge 0$, such that
for any subpath $p$ of $q$ the inequality $\lambda ||p|| - c \le d(p_-,p_+)$ holds.\\
In a hyperbolic space quasigeodesics and geodesics with same ends are mutually close~:

\begin{lemma} \label{close} {\normalfont (\cite[5.6,5.11]{Ghys},\cite[3.3]{Mihalik})}
There is a constant $\nu=\nu(\delta,\lambda,c)$ such that for any
$(\lambda,c)$-quasigeodesic path $p$ in $\Gamma(G,{\cal A})$ and a geodesic
$q$ with $p_- = q_-$, $p_+ = q_+$, one has $p \subset {\cal O}_\nu(q)$ and $q \subset
{\cal O}_\nu(p)$. \end{lemma}

An important property of cyclic subgroups in a hyperbolic group states

\begin{lemma} \label{power} {\normalfont (\cite[8.21]{Ghys},\cite[3.2]{Mihalik})}
For any word $w$ representing an element $g \in G$ of infinite order  there
exist constants $\lambda >0$, $c \ge 0$, such that any path with a label $w^m$ in the
Cayley graph of $G$ is $(\lambda,c)$-quasigeodesic for arbitrary integer $m$.
\end{lemma}

In particular, it follows from lemmas \ref{close} and \ref{power} that any cyclic subgroup of a
hyperbolic group is quasiconvex.

Recall that a group $H$ is called {\it elementary} if it has a
cyclic subgroup $\langle h \rangle$ of finite index.
It is known that every element $g \in G$ of infinite order is contained in a unique maximal
elementary subgroup $E(g)$ of $G$ (\cite{Gromov},\cite{Olsh2}), and
$$E(g) = \{x \in G~|~\exists~n\in \N ~\mbox{such that}~ xg^nx^{-1} = g^{\pm n}\}.$$

Let $W_1,W_2,\dots,W_l$ be words in $\cal A$ representing elements $g_1,g_2,\dots,g_l$
of infinite order, where $E(g_i) \neq E(g_j)$ for $i \neq j$. The following lemma
will be useful:

\begin{lemma} \label{quasigeodesic} {\normalfont (\cite[Lemma 2.3]{Olsh2})}
There exist
$\lambda = \lambda(W_1,W_2,\dots,W_l)>0$, $c=c(W_1,W_2,\dots,W_l) \ge 0$ and
$N = N(W_1,W_2,\dots,W_l)>0$ such that any path $p$ in the Cayley graph \ga
with label $W_{i_1}^{m_1}W_{i_2}^{m_2} \dots W_{i_s}^{m_s}$ is $(\lambda,c)$-quasigeodesic
if $i_k \neq i_{k+1}$ for $k=1,2,\dots,s-1$, and $|m_k|>N$ for $k=2,3,\dots,s-1$
(each $i_k$ belongs to $\{1,\dots,l\}$).
\end{lemma}

If $X_1,X_2,\dots, X_n$ are points in \ga, the notation $X_1X_2
\dots X_n$ will be used for the geodesic $n$-gon with vertices
$X_i$, $i=1,\dots,n$, and sides $[X_i,X_{i+1}]$, $i=1,2,\dots
,n-1$, $[X_n,X_0]$. $[X_1,X_2,\dots,X_n]$ will denote the broken line
with these vertices in the corresponding order.

\begin{lemma} \label{brokenlines1} {\normalfont (\cite[Lemma 21]{Olsh1})} Let $p = [X_0,X_1,\dots,X_n]$ be a broken
line in $\Gamma(G,{\cal A})$ such that
$||[X_{i-1},X_i]|| > C_1$ $\forall~i=1,\dots,n$, and $(X_{i-1}|X_{i+1})_{X_i} \le C_0$ $\forall~i=1,\dots,n-1$,
where $C_0 \ge 14 \delta$, $C_1 > 12(C_0 + \delta)$. Then $p$ is contained in the closed $2C_0$-neighborhood
${\cal O}_{2C_0}([X_0,X_n])$ of the geodesic segment $[X_0,X_n]$.
\end{lemma}

\begin{lemma} \label{brokenlines2} In the conditions of Lemma \ref{brokenlines1}, $\|[X_0,X_n]\| \ge \|p\|/2$.
\end{lemma}

\underline{Proof.} Induction on $n$. If $n=1$ the statement is trivial. So, assume $n>1$. By the induction hypothesis
$\|[X_0,X_{n-1}]\| \ge \|q\|/2$ where $q$ is the broken line $[X_0,X_1,\dots,X_{n-1}]$. It is shown in the proof of \cite[Lemma 21]{Olsh1} that
our conditions imply $(X_0|X_n)_{X_{n-1}} \le C_0 + \delta$, hence
$$\|[X_0,X_n]\| = \|[X_0,X_{n-1}]\| + \|[X_{n-1},X_n]\| - 2(X_0|X_n)_{X_{n-1}} \ge $$
$$ \ge \|q\|/2 + \|[X_{n-1},X_n]\|/2 + C_1/2 -2(C_0+\delta) \ge \|p\|/2~.$$
Q.e.d. $\square$

%


\begin{lemma} \label{intersection} {\normalfont (\cite[Prop. 3]{Short})} Let $G$ be a group generated by a finite set
 $\cal A$. Let $A,B$ be
subgroups of $G$ quasiconvex with respect to $\cal A$. Then $A \cap B$ is quasiconvex with respect to $\cal A$.
\end{lemma}

\begin{lemma} \label{double} {\normalfont (\cite[Prop. 1]{Arzh})}
 Let $G$ be a hyperbolic group and $H$ a quasiconvex
subgroup of $G$ of infinite index. Then the number of double cosets of $G$ modulo $H$ is infinite.
\end{lemma}

\begin{lemma} \label{long} {\normalfont (\cite[Lemma 10]{Arzh},\cite[Lemma 1.3]{GMR})}
For any integer $m \ge 1$ and numbers $\delta,\varepsilon,C \ge 0$, there exists $A=A(m,\delta,\varepsilon,C) \ge 0$
with the following property.

Let $G$ be a $\delta$-hyperbolic group with a generating set containing at most $m$ elements and $H$
a $\varepsilon$-quasiconvex  subgroup of $G$. Let $g_1,\dots,g_n,s$ be elements of $G$ such that

{\normalfont (i)} cosets $Hg_i$ and $Hg_j$ are different for $i \neq j$;

{\normalfont (ii)} $g_n$ is a shortest representative of the coset $Hg_n$;

{\normalfont (iii)} $|g_i|_G \le |g_n|_G$ for $1 \le i < n$;

{\normalfont (iv)} for $i \neq n$, all the products $g_ig_n^{-1}$ belong to the same double coset $HsH$

with $|s|_G \le C$. \\
Then $n \le A=A(m,\delta,\varepsilon,C)$.
\end{lemma}

\begin{lemma} \label{non-periodic} {\normalfont (\cite[8.3.36]{Ghys})}
Any infinite subgroup of a hyperbolic group contains an element of infinite order.
\end{lemma}

Let $H$ be a subgroup of $G$. Inheriting the terminology from \cite{GMR} we will say
that the elements $\{g_i~|~1 \le i \le n\}$ of $G$ are {\it essentially distinct} (relatively
to $H$) if $Hg_i \neq Hg_j$
for $i \neq j$. Conjugates $g_i^{-1}Hg_i$ of $H$ in this case are called {\it essentially
distinct conjugates}.

\vspace{.15cm}
{\bf \underline{Definition.}} (\cite[Def. 0.3]{GMR}) The {\it width} of an infinite subgroup $H$ in $G$ is $n$
if there exists a collection of essentially distinct conjugates of $H$ such that the
intersection of any two elements of the collection is infinite and $n$ is maximal possible.
The width of a finite subgroup is defined to be $0$.

\begin{lemma} \label{width} {\normalfont (\cite[Main Thm.]{GMR})}
A quasiconvex subgroup of a hyperbolic group has a finite width.
\end{lemma}

\begin{lemma} \label{quasiunion} {\normalfont (\cite[Cor. 1,Lemma 2.1]{hyp})}
In a hyperbolic group a finite union of quasiconvex products is a quasiconvex set.
\end{lemma}

\begin{lemma} \label{hypthm1} {\normalfont (\cite[Thm. 1]{hyp})}
Suppose $G_1, \dots, G_n$, $H_1, \dots, H_m$ are quasiconvex sub\-groups
of the group $G$, $n,m \in \N$; $f,e \in G$.
Then there exist numbers \\ $r, t_1,\dots,t_r \in \N \cup \{0\}$ and $f_l,\alpha_{lk},
\beta_{lk} \in G$, $k =1,2,\dots, t_l$ (for every
fixed $l$),\\ $l =1,2,\dots, r$, such that
$$fG_1 G_2 \cdot \dots \cdot G_n \cap eH_1 H_2
\cdot \dots \cdot H_m = \bigcup_{l=1}^r f_l S_l $$
where for each $l$, $t = t_l$, there are indices $1 \le i_1\le i_2\le \dots \le i_t \le n$,
$1 \le j_1\le $ $ \le j_2\le \dots \le j_t \le m$~:
$$ \quad S_l = (G_{i_1}^{\alpha_{l1}} \cap H_{j_1}^{\beta_{l1}}) \cdot
\dots \cdot (G_{i_t}^{\alpha_{lt}} \cap H_{j_t}^{\beta_{lt}}). $$
\end{lemma}

\remark  \label{rem3} Observe that arbitrary quasiconvex product
$f_0G_1f_1G_2 \cdot \dots \cdot G_nf_n$
is equal to a "transformed" product $f G'_1 G'_2 \cdot \dots \cdot G'_n$ (which appears in the
formulation of Lemma \ref{hypthm1}) where
$G'_i = (f_{i} \cdot \dots \cdot f_{n})^{-1}G_i(f_{i} \cdot \dots \cdot f_{n})$,
$i=1,\dots,n$, are quasiconvex subgroups of $G$ by Remark \ref{rem2} and
$f = f_0f_1 \cdot \dots f_{n} \in G$.

\section{\large Properties of Quasiconvex Subgroups of Infinite Index}

\begin{lemma} \label{quadrangle}  Consider  a geodesic quadrangle $X_1X_2X_3X_4$  in \ga with \\
$d(X_2,X_3)>d(X_1,X_2)+d(X_3,X_4)$. Then there are points $U,V \in [X_2,X_3]$ such that
$d(X_2,U) \le d(X_1,X_2)$, $d(V,X_3) \le d(X_3,X_4)$ and  the geodesic subsegment
$[U,V]$ of $[X_2,X_3]$ lies $2\delta$-close to the side $[X_1,X_4]$.
\end{lemma}

\underline{Proof.} Since $(X_1|X_3)_{X_2} \le d(X_1,X_2)$ and $(X_1|X_4)_{X_3} \le d(X_3,X_4)$, one can choose
$U,V \in [X_2,X_3]$ satisfying $d(X_2,U) = (X_1|X_3)_{X_2}$, $d(X_3,V) = (X_1|X_4)_{X_3}$. The triangle
$X_1X_3X_2$ is $\delta$-thin, therefore, after taking $V'\in [X_1,X_3]$ at distance $d(X_3,V)$ from
$X_3$, one obtains $[U,V] \subset \O_{\delta}([X_1,V'])$. Finally, since $V'$ is the special point
of triangle $X_1X_3X_4$ by construction, $[X_1,V']$ is in the closed $\delta$-neighborhood of the side $[X_1,X_4]$, and
thus, $[U,V] \subset \O_{2\delta}([X_1,X_4])$. $\square$

\begin{lemma} \label{conj} Let $A$ be an infinite $\varepsilon$-quasiconvex set in $G$  and $g \in G$.
Then if the intersection $A \cap gAg^{-1}$ is infinite, there exists an element $r\in G$ with
$|r|_G \le 4\delta+ 2\varepsilon + 2\varkappa$ such that $g \in ArA^{-1}$, where $\varkappa$ is the length
of a shortest element from $A$.
\end{lemma}

\underline{Proof.} Note, at first, that for every $a \in A$ the geodesic segment $[1_G,a]$  belongs
to a closed $(\delta+\varepsilon+\varkappa)$-neighborhood of $A$ in \ga. Indeed, pick $b \in A$ with
$d(1_G,b)=|b|_G = \varkappa$ and consider the geodesic triangle $1_Gab$. Using $\delta$-hyperbolicity of
the Cayley graph one achieves
$$[1_G,a] \subset \O_\delta([a,b] \cup [1_G,b]) \subset \O_{\delta+\varkappa}([a,b]) \subset
\O_{\delta+\varkappa+\varepsilon}(A)~.$$

\begin{figure}
\input{pic1.tex}
\begin{center}
Figure 1
\end{center}

\end{figure}

By the conditions of the lemma there is an element $a_1 \in A$ such that $ga_1g^{-1} = a_2 \in A$ and
$|a_1|_G>2|g|_G$. Set $X_1 =1_G$, $X_2 = g$, $X_3=ga_1$,$X_4=a_2$ (Fig. 1). Then $d(X_2,X_3) = |a_1|_G$,
$d(X_1,X_2)=|g|_G=|a_2^{-1}ga_1|_G=d(X_3,X_4)$ and in the geodesic quadrangle
$X_1X_2X_3X_4$ one has $d(X_2,X_3)>d(X_1,X_2)+d(X_3,X_4)$ and, so, by lemma 1 there exist $x \in [X_1,X_4]$,
$y \in [X_2,X_3]$ with $d(x,y) \le 2\delta$. As we showed above, $[1_G,a_i] \subset \O_{\delta+\varkappa+\varepsilon} (A)$
for $i=1,2$, hence there is $\alpha \in A$ such that
$d(\alpha,x) \le \delta+\varepsilon+\varkappa$. A left shift is an isometry of \ga, thus,
$[X_2,X_3]=[g,ga_1] \subset \O_{\delta+\varkappa+\varepsilon} (gA)$ and we can obtain  an
element $\beta \in A$ such that $d(y,g\beta) \le \delta+\varepsilon+\varkappa$.

Consider the broken line $q = [X_1,\alpha,g\beta,g]$ in \ga; then $elem(q) = g$ in $G$.
$d(\alpha,g\beta) \le 4\delta+ 2\varepsilon + 2\varkappa$
by construction, hence we achieved $g=elem(q) = \alpha \cdot r \cdot \beta^{-1}$ where $r = elem([\alpha,g \beta])$,
$|r|_G \le d(\alpha,g\beta) \le 4\delta+ 2\varepsilon + 2\varkappa$. $\square$

\vspace{.15cm}
For the case when $A$ is a quasiconvex subgroup of a hyperbolic group $G$,
Lemma \ref{conj} was proved in \cite[Lemma 1.2]{GMR}.

\begin{cor} \label{cor2}
Suppose $H$ is a quasiconvex subgroup of infinite index in a hyperbolic group $G$.
Then  $H$ contains no infinite normal subgroups of $G$.
\end{cor}

\underline{Proof.} Indeed, assume $N \unlhd G$ and $N \subset H$.
By Lemma \ref{double}, there is a double coset $HrH$, $r \in G$, with the length of a shortest
representative greater than $(4\delta+ 2\varepsilon)$ ($\varepsilon$ is the quasiconvexity constant of $H$).
Thus, according to the Lemma \ref{conj}, $N \subset H \cap rHr^{-1} $ is finite. $\square$

\vspace{.15cm}

\begin{lemma} \label{newdouble} Let $G$ be a $\delta$-hyperbolic group, $H$ and $K$ -- its subgroups
with $H$  quasiconvex. If ~~$\displaystyle K \subset \bigcup_{j=1}^N Hs_jH$~~ for some $s_1,\dots,s_N \in G$ then $K \preceq H$,
i.e. $|K:(K \cap H)| < \infty$ .
\end{lemma}

\underline{Proof.} By contradiction, assume $K = \bigsqcup_{i=1}^\infty (K \cap H)x_i$ -- disjoint union of
right cosets with $x_i \in K$ for all $i \in \N$. For every $i$ choose  a shortest representative $g_i$ of
the coset $Hx_i$ in $G$. Then for arbitrary  $i \neq k$, $Hg_i=Hx_i \neq Hx_k=Hg_k$ and
$x_ix_k^{-1} \in Hs_jH$ for some $j \in \{1,2,\dots,N\}$, $j=j(i,k)$, hence $g_ig_k^{-1} \in Hs_jH$.

Let $A_j$ be the constants corresponding to $Hs_jH$, $j=1,\dots,N$, from Lemma \ref{long}. Pick a natural
number $n > \sum_{j=1}^N A_j$ and consider $g_1,g_2,\dots,g_n$. Without loss of generality, assume
$|g_n|_G \ge |g_i|_G$ for $1 \le i < n$.

By the choice of $n$, there exits $l \in \{1,\dots,N \}$ such that
$$card\left\{i \in \{1,2,\dots,n-1\}~|~ g_ig_n^{-1} \in Hs_lH \right\} \ge A_l~.$$
 This leads to a contradiction
to Lemma \ref{long}. Q.e.d. $\square$

\begin{lemma} \label{infinite} Suppose $H$ is a quasiconvex subgroup of a hyperbolic group $G$ and $K$
is a subgroup of $G$ with $|K:(K \cap H)| = \infty$. Then there is  an element
$x \in K$ of infinite order such that the intersection $\langle x \rangle_\infty \cap H$ is trivial.
\end{lemma}

\underline{Proof.} Observe that our conditions imply that $K$ is infinite.
If $K \cap H$ is finite, the statement follows from Lemma \ref{non-periodic}
applied to $K$.

So, suppose $K \cap H$ is infinite; hence there is an element $y \in K \cap H$ of infinite order.
$|K:(K \cap H)| = \infty$, therefore, applying lemmas \ref{newdouble} and \ref{conj} one obtains
$g \in K$ such that $H \cap gHg^{-1}$ is a finite subgroup of $H$.
Thus $x = gyg^{-1} \in gHg^{-1}$ satisfies the needed property
($x$ is an element of $K$ because $g,y \in K$). $\square$

\begin{lemma} \label{quasintersection} Let $G$ be a $\delta$-hyperbolic group with respect to some finite
generating set $\cal A$, and let $H_i$ be $\varepsilon_i$-quasiconvex subgroups of $G$, $i=1,2$.
If one has $\sup \{(h_1|h_2)_{1_G}~:~h_1 \in H_1,h_2 \in H_2 \} = \infty$ then
$card(H_1\cap H_2)=\infty$.
\end{lemma}

\underline{Proof.} Define a finite subset $\Delta$ of $G$ by
$\Delta= \{g \in H_1H_2~:~|g|_G \le \delta+\varepsilon_1+\varepsilon_2\}$. For each
$g \in \Delta$ pick a pair $(x,y) \in H_1 \times H_2$ such that $x^{-1}y =g$, and
let $\Omega \subset H_1 \times H_2$ denote the (finite) set of the chosen pairs.
Define $\Omega_1$
to be the projection of $\Omega$ on $H_1$, i.e. $\Omega_1 = \{x \in H_1~|~\exists~ y \in H_2: (x,y) \in \Omega\}$.

By construction, $card(\Omega_1) < \infty$, and thus, $D \stackrel{def}{=}
max\{|h|_G~:~h \in \Omega_1 \} <\infty$.

\begin{figure}
\begin{center}
\input{pic2.tex}
\end{center}

\begin{center}
Figure 2
\end{center}

\end{figure}

Assume $h_1 \in H_1$ and $h_2 \in H_2$ and consider the geodesic triangle $1_Gh_1h_2$ in
\ga (Figure 2). Let $P,Q$ be its "special" points on the sides $[1_G,h_1]$ and $[1_G,h_2]$ correspondingly.
Since the triangles in \ga are $\delta$-thin and $H_i$ are $\varepsilon_i$-quasiconvex, $i=1,2$,
we have $d(P,Q) \le \delta$, and there exist $\h_i \in H_i$, $i=1,2$, with
$d(\h_1,P) \le \varepsilon_1$, $d(\h_2,Q) \le \varepsilon_2$. By the triangle inequality
$d(\h_1,\h_2) = |\h_1^{-1}\h_2|_G \le \delta+\varepsilon_1+\varepsilon_2$, therefore, there
is a pair $(x,y) \in \Omega$ such that $\h_1^{-1}\h_2 = x^{-1}y$, and so
\begin{equation} \label{h1} \h_1x^{-1} = \h_2 y^{-1} \in H_1 \cap H_2~.
\end{equation}

From our construction it also follows that
\begin{equation} \label{h2}
|\h_1x^{-1}|_G \ge |\h_1|_G-|x|_G \ge d(1_G,P)-\varepsilon_1 - D = (h_1|h_2)_{1_G} -
\varepsilon_1-D~ .
\end{equation}

Hence, if $\sup \{(h_1|h_2)_{1_G}~:~h_1 \in H_1,h_2 \in H_2 \} = \infty$ then, as we see from
(\ref{h1}) and (\ref{h2}), $H_1 \cap H_2$ has elements of arbitrary large lengths, and thus, it
is infinite. $\square$

\vspace{.15cm}
Suppose $H_1,H_2,\dots,H_s$ are quasiconvex subgroups of the group $G$, and $K$ is a subgroup of $G$
satisfying $|K:(K \cap H_j)| = \infty$, $j=1,2,\dots,s$.
\begin{lemma} \label{notinunion} There exists an element of infinite order $x \in K$
with the property $\langle x \rangle_\infty \cap \left( H_1 \cup H_2 \cup \dots \cup H_s \right)
= \{1_G\}$.
\end{lemma}

\underline{Proof.} By Lemma \ref{infinite} for every $i=1,2,\dots,s$ there is an
element $\tilde x_i \in K$ of infinite order such that $\tilde x_i^n \notin H_i$ for every
$n \in \Z \backslash \{0\}$. If for some $i,j \in \{1,\dots,s\}$, $i \neq j$, one has
$E(\tilde x_i) = E(\tilde x_j)$, then $\langle \tilde x_i \rangle \cap \langle \tilde x_j \rangle$
is non-trivial, therefore one can remove $\tilde x_j$ from the collection
$\{\tilde x_1,\tilde x_2,\dots,\tilde x_s\}$
because in this case one has $\langle \tilde x_i \rangle_\infty \cap H_j = \{1_G\}$ (since any two non-trivial
subgroups of $\Z$ intersect non-trivially).
Thus, after performing this procedure a finite number of times, we will obtain a collection
of elements of infinite order $\{x_1,x_2,\dots,x_r\} \subset G$, $r \le s$, satisfying the properties:
$E(x_i) \neq E(x_j)$ for $i \neq j$, and for every $k =1,2,\dots,s$ there is
$i_k \in\{1,\dots,r\}$ such that $\langle x_{i_k} \rangle_\infty \cap H_k = \{1_G\}$.

Now, for each $i =1,2,\dots,r$ and $k =1,2,\dots,s$, define
$$\alpha_{ik} = \left\{ \begin{array}{lc} min \{m \in \N~|~x_i^m \in H_k\},
~&\mbox{if}~\langle x_i \rangle \cap H_k \neq \{1_G\} \\
0,~ &\mbox{otherwise}
\end{array} \right. .$$
Denote
$$ \alpha_i =\left\{ \begin{array}{lc} 1,
~&\mbox{if}~\alpha_{ik}=0 ~\mbox{for every $k=1,\dots,s$} \\
l.c.m.\{\alpha_{ik}~|~1\le k \le s,\alpha_{ik} >0 \},~ &\mbox{otherwise}
\end{array} \right. .$$

$\alpha_i \in \N$, so one can set $y_i = x_i^{\alpha_i}$, $i=1,2,\dots,r$.
Then for any distinct $i,j \in \{1,\dots,r\}$, and any $k =1,\dots,s$,
$E(x_i)=E(y_i) \neq E(y_j)=E(x_j)$,
and

\vspace{.15cm}
{\bf (i)} either $y_i \in H_k$ or $\langle y_i \rangle_\infty \cap H_k = \{1_G\}$ (by
construction, the latter holds for $i=i_k$).

\vspace{.15cm}
At last, for every natural $n$ we define $z_n = y_1^ny_2^n \cdot\dots\cdot y_r^n \in K$.

Assume, by contradiction, that for each $n \in N$ there exits $l=l(n) \in \N$ such that
$z_n^l \in H_1\cup \dots \cup H_s$ (this obviously holds if $z_n$ has a finite order in $G$).
Then there is an index $k_0 \in \{1,2,\dots,s\}$ such that $(z_n)^{l(n)} \in H_{k_0}$ for
infinitely many $n \in \N$. Without loss of generality, assume  $k_0 =1$.

Let $w_1,w_2,\dots,w_r$ be words over the alphabet $\cal A$ representing the elements
$y_1,y_2,\dots,y_r$ correspondingly.
We apply Lemma \ref{quasigeodesic} to obtain $\lambda >0$,
$c \ge 0$ and $N>0$ (depending on $w_1,w_2,\dots,w_r$) such that any path
$p$ in the Cayley graph \ga with label
$(w_{1}^{n}w_{2}^{n} \dots w_{r}^{n})^{l(n)}$ is $(\lambda,c)$-quasigeodesic
if $n>N$. Let $\nu = \nu(\delta,\lambda,c)$ be the constant from Lemma \ref{close}, and
$\varepsilon$ be the quasiconvexity constant of $H_1$.

Let $p_n$ be the path in \ga starting at $1_G$ labelled by
$(w_{1}^{n}w_{2}^{n} \dots w_{r}^{n})^{l(n)}$, $n \in \N$, $n >N$. $p_n$ ends at the element
$z_n^{l(n)}$ which belongs to $H_1$ for infinitely many $n$. Then
$p_n \subset \O_{\nu+\varepsilon}(H_1)$ for infinitely many $n$.

Denote $t = min\{i~|~1\le i \le r, y_i \notin H_1 \}$ (such $t$ exists by construction of
$y_i$). By definition, $y_1^{n}y_2^n \cdot \dots \cdot y_t^n$ lies on $p$, therefore,
$d(y_1^{n}y_2^n \cdot \dots \cdot y_t^n,H_1) \le \nu+\varepsilon$ for infinitely many $n$,
i.e. for those $n$ one can find elements $a_n \in G$ such that
$y_1^{n}y_2^n \cdot \dots \cdot y_t^na_n \in H_1$ and $|a_n|_G \le \nu+\varepsilon$.
Remark that $y_1,\dots,y_{t-1} \in H_1$ by definition, hence $y_t^n a_n = h_n \in H_1$
for infinitely many $n \in \N$.

Finally, since one has
$$(y_t^n|h_n)_{1_G} = 1/2(|y_t^n|_G+|h_n|_G-|a_n|_G)\ge
1/2(|y_t^n|_G-\nu-\varepsilon)$$ for infinitely many $n$, we apply Lemma \ref{power} to achieve
$\sup\{(y_t^n|h_n)_{1_G}\} = \infty$. Hence, by Lemma \ref{quasintersection},
$card(\langle y_t \rangle_\infty \cap H_1) = \infty$, thus, from $(i)$ it follows that
$y_t \in H_1$ -- a contradiction.

Therefore, there exist $n \in \N$ such that $z_n^m \notin H_1 \cup \dots \cup H_s $ for
every $m \in N$, consequently, $x=z_n \in K$ has an infinite order and
$\langle x \rangle_\infty \cap (H_1 \cup \dots \cup H_s) = \{1_G\}$.
The proof of the lemma is finished. $\square$

\begin{prop} \label{thm1}
Suppose $H$ is a quasiconvex subgroup of a hyperbolic group $G$ and $K$ is any subgroup of
$G$ that satisfies $|K:(K \cap H^g)| = \infty$ for all $g\in G$. Then there
exists an element $x \in K$ having infinite order, such that
$\langle x \rangle_\infty \cap H^G = \{1_G\}$.
\end{prop}

\underline{Proof} of Proposition \ref{thm1}.
Observe that the subgroup $K$ is infinite because $|K:(K \cap H)| = \infty$, thus
applying Lemma \ref{non-periodic} we obtain an element $h \in K$ of infinite order.

Assume, by contradiction, that for every $x \in K$ there exist $l=l(x) \in \N$ and $g=g(x) \in G$
such that $gx^lg^{-1} \in H$. In particular, $h^{l_0} \in g_0^{-1}Hg_0$ for some $l_0 \in \N$,
$g_0 \in G$.

Take an arbitrary $y \in K$ of infinite order. If $E(y)=E(h)$ then
$y^m \in \langle h \rangle$ for some $m \in \N$, hence $y^{ml_0} \in g_0^{-1}Hg_0$.

\begin{figure}
\input{pic3.tex}

\begin{center}
Figure 3
\end{center}

\end{figure}

If $E(y) \neq E(h)$, choose words $w_1, w_2$ over the alphabet $\cal A$ representing $y$ and
$h$. Then by Lemma \ref{quasigeodesic} there exist
$\lambda = \lambda(w_1,w_2)>0$, $c=c(w_1,w_2) \ge 0$ and
$N = N(w_1,w_2)>0$ such that any path $p$ in the Cayley graph \ga
with label $(w_1^nw_2^n)^l$ is $(\lambda,c)$-quasigeodesic
if $n>N$, for every $l \in \N$.

Let $\nu=\nu(\delta,\lambda,c)$ be the constant from Lemma \ref{close}  and let $\varepsilon$
denote the quasiconvexity constant for $H$.

By our assumption for every $n \in \N$, $n>N$, there are $l_n \in \N$ and $g_n \in G$
satisfying $g_n (y^nh^n)^{l_n} g_n^{-1} \in H$, then $g_n (y^nh^n)^{kl_n} g_n^{-1} \in H$
$\forall~k \in \N$. Consider a path $p_n$ in \ga with
$(p_n)_{-} = g_n$ and $lab(p_n) \equiv (w_1^nw_2^n)^{kl_n}$ and a geodesic quadrangle
$X_1X_2X_3X_4$ in \ga where $X_1 = 1_G$, $X_2 = g_n$, $X_3 = g_n(y^nh^n)^{kl_n}$,
$X_4 = g_n (y^nh^n)^{kl_n} g_n^{-1} \in H$.

Now, if $n>N$, the path $p_n$ is $(\lambda,c)$-quasigeodesic and $\nu$-close to $[X_2,X_3]$,
therefore, applying Lemma \ref{quadrangle}, we conclude that for sufficiently large $k$
(compared to $\nu$ and $|g_n|_G$) there is a subpath $q_n$ of $p_n$ labelled by
$w_1^nw_2^n$ which lies in  the closed $C=C(\delta,\nu)$-neighborhood of $[X_1,X_4]$, consequently,
$q_n \subset \O_{C+\varepsilon}(H)$. Note that $C$ depends only on $\delta,y,h$ but does not
depend on $n,k$ and $g_n$.

Hence, for each natural $n>N$ we obtained elements
$u_n,v_n,z_n \in G$ with $|u_n|_G,|v_n|_G,|z_n|_G \le C+\varepsilon$ and $u_ny^nv_n^{-1} \in H$,
$v_nh^nz_n^{-1} \in H$ (see Figure 3). There are infinitely many such $n$
and only finitely many possible $u_n,v_n,z_n$, hence for some indices
$i,j \in \N$, $i<j$, one will have $u_i=u_j$, $v_i=v_j$, $z_i=z_j$. Thus,
$u_iy^iv_i^{-1} \in H$, $u_iy^jv_i^{-1} \in H$, implying $v_iy^{j-i}v_i^{-1} \in H$.
Similarly, $v_ih^{j-i}v_i^{-1} \in H$.

Thus for each $y \in K$ of infinite order we found an element $a =a(y)\in G$ and $l = l(y) \in \N$
satisfying $ay^la^{-1} \in H$ and $ah^la^{-1} \in H$
(in the case when $E(y)=E(h)$, $a=g_0$, $l=ml_0$).

This implies that for arbitrary $y_1,y_2 \in K$ with $o(y_i)=\infty$, $i=1,2$, we have
$a_i =a(y_i)\in G$ and $t_i = t(y_i) \in \N$ such that $y_i^{t_i} \in a_i^{-1}Ha_i$ and
$h^{t_i} \in a_i^{-1}Ha_i$, $i=1,2$. Therefore, $h^{t_1t_2} \in a_1^{-1}Ha_1 \cap a_2^{-1}Ha_2$, and, thus,
this intersection is infinite (because this subgroup contains an element of infinite order).
But by Lemma \ref{width} there can be only finitely many of such conjugates of $H$ in $G$,
therefore, there are $a_1,a_2,\dots,a_s \in G$ such that
for every element $y \in G$ of infinite order one has
$$y^l \in \bigcup_{i=1}^s a_i^{-1}Ha_i$$ for some $l \in \N$,
which contradicts to Lemma \ref{notinunion} because of Remark \ref{rem2} and the conditions of the proposition.

Hence there is an element $x \in K$ such that
for any $l \in \N$, $x^l \notin H^G$ . It follows that $x$ has infinite order and
$\langle x \rangle_\infty \cap H^G = \{1_G\}$. $\square$

\section{\large Proofs of Theorems \ref{prop1},\ref{prop2}}

\underline{Proof of Theorem\ref{prop1}.} The direction (b) $\Rightarrow$ (a) is trivial, so let's
focus on the direction (a) $\Rightarrow$ (b).

Using Proposition \ref{thm1} for every $i=1,2,\dots,s$ one finds
$\tilde x_i \in K$ such that $o(\tilde x_i) = \infty$ and $\langle \tilde x_i \rangle_\infty
\cap H_i^G = \{1_G\}$.
Now as in the proof of the Lemma \ref{notinunion} we narrow down this collection to
$\{ x_1,x_2,\dots,x_r \}$ satisfying the properties: $E(x_i) \neq E(x_j)$ for $i \neq j$,
and for every $k=1,2,\dots,s$ there is $i_k \in \{1,2\dots,r\}$ such that
$\langle  x_{i_k} \rangle_\infty \cap H_k^G = \{1_G\}$.

Assuming the contrary of the statement, for every natural number $n$ and
$z_n = x_1^nx_2^n \cdot \dots \cdot x_r^n \in K$ we obtain $l_n \in \N$ such that
$z_n^{l_n} \in H_1^G \cup \dots \cup H_s^G$. Thus, there is an index $k_0 \in \{1,2,\dots,s\}$
such that $z_n^{l_n} \in H_{k_0}^G$ for every $n \in \Delta$ where $\Delta$ is an infinite
subset of $\N$. Without loss of generality, assume $k_0 =1$.

Thus, for every $n \in \Delta$ there exists $g_n \in G$ such that
$g_nz_n^{l_n}g_n^{-1} \in H_1$, hence $g_nz_n^{kl_n}g_n^{-1} \in H_1$ for any $k \in \N$.
Now, as in the proof of the Proposition \ref{thm1}, we take words $w_1,w_2,\dots,w_r$ representing
$x_1,x_2,\dots,x_r$ and apply Lemma \ref{quasigeodesic} to the path $p_n$ in \ga starting at
$g_n$ and labelled by $(w_1^nw_2^n \dots w_r^n)^{kl_n}$. Hence, for every sufficiently large
$n \in \Delta$ and $k \in \N$ we find a subpath $q_n$ of $p_n$ labelled by $w_1^n \dots w_r^n$
which is $(C+\varepsilon)$-close to $H_1$ where $C=C(\delta,w_1,\dots,w_r)$ (but is independent
of $n$,$k$ and $g_n$) and $\varepsilon$ is the quasiconvexity constant of $H_1$.

Again, similarly to the proof of Proposition \ref{thm1}, for each $i=1,2,\dots,r$ we obtain
an element $v_i \in G$ and $m_i \in \N$ such that $v_ix_i^{m_i}v_i^{-1} \in H_1$. In
particular that should hold for $i=i_1$. The contradiction achieved finishes the proof
of the theorem. $\square$

\vspace{.15cm}
\underline{Proof of Theorem\ref{prop2}.} Assume the contrary. Then
using Theorem\ref{prop1} we obtain an element  $x \in K$, $o(x)=\infty$,
such that $\langle x \rangle_\infty \cap (H_1^G \cup H_2^G \cup \dots \cup H_s^G)$ is trivial,
hence $\langle x \rangle^G \cap H_j^G = \{1_G\}$ for every $j=1,2,\dots,s$.
%

By the conditions $U = \bigcup_{k=1}^q P_k$ where each $P_k$ is a quasiconvex product, $k=1,\dots,q$.
Since any cyclic subgroup is quasiconvex in $G$ (by lemmas \ref{power},\ref{close}), and
in view of Remark \ref{rem3}, application of the Lemma \ref{hypthm1} to the intersection
$$\langle x \rangle=\langle x \rangle \cap U = \bigcup_{k=1}^q \bigl( \langle x \rangle \cap P_k \bigr)$$
shows that it is finite (because each product $S_l$ will be trivial). A contradiction
with $o(x)=\infty$. Hence, the theorem is proved. $\square$

%
%


\begin{cor} \label{non-quasidense} Let $U$ be a
finite union of quasiconvex products having infinite index in a
hyperbolic group $G$. Then $U$ is not quasidense but its
complement $U^{(c)} = G \backslash U$ is quasidense in $G$.
\end{cor}

\underline{Proof.} Let $g_1,g_2, \dots,g_n \in G$.
It is easy see that the sets $\bigcup_{i=1}^n Ug_i$ and $U^{-1}U$ are
finite unions of quasiconvex products of infinite index.

Thus, the fact that $U \subset G$ is not quasidense follows from the definition of a quasidense subset and
Theorem\ref{prop2} applied to the case when $K=G$.

By Theorem\ref{prop2}, there exists $y \in G$ such that $y \notin U^{-1}U$. 
Consider an arbitrary $x \in G$. If $x \in U$ then $xy \in U^{(c)}$ (because of the choice of $y$) hence
$$ G \subseteq U^{(c)} \cup U^{(c)}y^{-1}~. ~~\square$$

\section{\large Boundary and Limit Sets}
Let $X$ be a proper geodesic metric space with metric $d(\cdot,\cdot)$. Assume also that
$X$ is $\delta$-hyperbolic for some $\delta \ge 0$.
Further in this paper we will need the construction of Gromov {\it boundary} $\partial X$ for the space $X$
(for more detailed theory the reader is referred to the corresponding chapters in \cite{Ghys},\cite{Bridson}).
The points of the boundary are equivalence
classes of geodesic rays $r: [0,\infty) \to X$ where rays $r_1,r_2$ are equivalent if
$\sup \{d(r_1(t),r_2(t))\} < \infty$ ($\Leftrightarrow~h(r_1,r_2) < \infty$ -- Hausdorff distance between the
images of these rays).

For another definition of the boundary, fix a basepoint $p \in X$. A sequence $(x_i)_{i\in \N} \subset X$ is called
{\it converging to infinity} if
$$\lim_{i,j \to \infty} (x_i|x_j)_p = \infty~.$$

Two sequences $(x_i)_{i\in \N},(y_j)_{i\in \N}$ converging to infinity are said to be equivalent
if $$\lim_{i \to \infty} (x_i|y_i)_p = \infty~.$$
The points of the boundary $\partial X$ are identified with the equivalence classes of sequences
converging to infinity. It is easy to see that this definition does not depend on the choice
of a basepoint. If $\alpha$ is the equivalence class of $(x_i)_{i\in \N}$ we will write
$\displaystyle \lim_{i \to \infty} x_i = \alpha$.

It is known that the two objects defined above are homeomorphic
through the map sending a geodesic ray $r:[0,\infty) \to X$ into the sequence $\bigl(r(i)\bigr)_{i \in \N}$.

For any two distinct points $\alpha, \beta \in \partial X$ there exists at least one bi-infinite
geodesic $r:(-\infty,+\infty) \to X$ such that $\displaystyle \lim_{i \to \infty}r(-i) = \alpha$ and
$\displaystyle \lim_{i \to \infty}r(i) = \beta$. We will say that this geodesic joins $\alpha$
and $\beta$; it will be denoted $(\alpha,\beta)$.

The spaces $\partial X$ and $X \cup \partial X$ can be topologized so that they become compact
and Hausdorff (see \cite{Mihalik},\cite{Ghys}).

Every isometry $\psi$ of the space $X$ induces a homeomorphism of $\partial X$ in a natural way:
for every equivalence class of geodesic rays $[r] \in \partial X$ choose a representative $r: [0,\infty) \to X$
and set $\psi([r])=[\psi \circ r]$.

For a subset $A \subseteq X$ the {\it limit set} $\Lambda (A)$ of $A$ is the collection of the points
$\alpha \in \partial X$ that are limits of the sequences from $A$.

Let $\Omega$ be a subset of $\partial X$ containing at least two distinct points. We define
the {\it convex hull} $CH(\Omega)$ of $\Omega$ to be the set of all points in $X$ lying on bi-infinite geodesics
that join elements from $\Omega$.


Below we list some known properties of limit sets and convex hulls:
\begin{lemma} \label{convexhulls} {\normalfont (\cite[Lemmas 3.2,3.6]{K-S})}
Let $\Omega$ be an arbitrary subset of $\partial X$ having at least two elements. Then

(a) $CH(\Omega)$ is $\varepsilon$-quasiconvex where $\varepsilon \ge 0$ depends only on $\delta$;

(b) If the subset $\Omega$ is closed then $\Lambda\left(CH(\Omega)\right) = \Omega.$

\end{lemma}

In this paper our main interest concerns hyperbolic
groups, so further we will assume that the space $X$ is the Cayley
graph \ga of some $\delta$-hyperbolic group $G$ with a fixed
finite generating  set $\cal A$. Because of the natural embedding
of $G$ (as a metric subspace) into \ga, we will identify subsets
of $G$ with subsets of its Cayley graph. The boundary $\partial
G$, by definition, coincides with the boundary of \ga.

Left multiplication by elements of the group induces the isometric action of $G$ on \ga. Hence,
$G$ acts homeomorphically on the boundary $\partial G$ as described above.

If $g \in G$ is an element of infinite order in $G$ then the sequences $(g^i)_{i\in \N}$ and
$(g^{-i})_{i\in \N}$ converge to infinity and we will use the notation
$$\lim_{i \to \infty} g^i = g^\infty,~\lim_{i \to \infty} g^{-i} = g^{-\infty}~.$$

\begin{lemma} \label{limitsets} {\normalfont (\cite{K-S},\cite{Swenson})}
Suppose $A,B$ are arbitrary subsets of $G$, $g \in G$. Then

(a) $\Lambda (A) = \emptyset$ if and only if $A$ is finite;

(b) $\Lambda (A)$ is a closed subset of the boundary $\partial G$;

(c) $\Lambda(A \cup B) = \Lambda (A) \cup \Lambda (B)$;

(d) $\Lambda(Ag) = \Lambda(A)$, $g \circ \Lambda(A)=\Lambda(gA)$;

(e) If $A \preceq B$ then $\Lambda(A) \subseteq \Lambda(B)$. Hence, $A \approx B$ implies $\Lambda (A) = \Lambda (B)$.
\end{lemma}

\noindent (b),(c) and (d) are easy consequences of the definition and (a) is obtained after a standard
application of Ascoli theorem; (e) follows from (c) and (d).

If $H$ is a subgroup of $G$, it is known that $\Lambda H$ is either empty (if $H$ is finite) or
consists of two distinct points (if $H$ is infinite elementary) or is uncountable (if $H$ is
non-elementary). In the second case, i.e. when there exists $g \in H$ such that $o(g) = \infty$ and
$|H:\langle g \rangle| < \infty$, $\Lambda H = \{g^{\infty},g^{-\infty} \}$.

\begin{lemma} \label{invariant} {\normalfont (\cite[Lemma 3.3]{K-S})} If $H$ is an infinite subgroup of $G$ then
$\Lambda(H)$ contains at least
two distinct points and the sets $\Lambda(H)$, $CH\bigl(\Lambda(H)\bigr)$ are $H$-invariant, i.e.
for every $h \in H$ ~~$h\circ \Lambda(H)=\Lambda(H)$, $h\cdot CH\bigl(\Lambda(H)\bigr) = CH\bigl(\Lambda(H)\bigr)$.

\end{lemma}

As the hyperbolic group $G$ acts on its boundary, for every subset $\Omega \subset \partial G$ one can
define the stabilizer subgroup by $St_G(\Omega) = \{g \in G~|~g \circ \Omega = \Omega\}$.
For our convenience, we set $St_G(\emptyset) = G$.

It is proved in \cite[thm. 8.3.30]{Ghys} that for any point $\alpha \in \partial G$~~ $St_G(\{\alpha\})$ is
an elementary subgroup of the group $G$ (in fact, if $\alpha=g^\infty$ for some element
of infinite order $g \in G$ then $$St_G(\{\alpha\})=E^+(g)\stackrel{def}{=}
\{x \in G~|~\exists~ n \in \N ~\mbox{such that}~xg^nx^{-1}=g^n\} \le E(g);$$
otherwise the subgroup $St_G(\{\alpha\})$ is finite). In addition, if $g \in G$, $o(g)=\infty$, then
$St_G(\{g^\infty,g^{-\infty}\}) = E(g)$.

\vspace{.15cm}
\remark  \label{rem4} For an arbitrary subset $A$ of $G$ $Comm_G(A) \subseteq St_G\bigl(\Lambda(A)\bigr)$.

\vspace{.15cm}
Indeed, if $g \in Comm_G(A)$, then $gA \approx A$, hence after applying claims (d),(e) of Lemma \ref{limitsets},
we obtain $g \circ \Lambda(A)=\Lambda(gA)=\Lambda(A)$, i.e. $g \in St_G\bigl(\Lambda(A)\bigr)$.

\vspace{.15cm}
\remark \label{LCH-sets} Suppose $\Omega \subseteq \partial G$ has at least two distinct points. Denote by
$cl(\Omega) \subseteq \partial G$ the closure of $\Omega$ in the topology of the group boundary. Then
$\Lambda\left(CH(\Omega)\right) = cl(\Omega).$

\vspace{.15cm}
Indeed, since $CH(\Omega) \subseteq CH\left(cl(\Omega)\right)$ we obtain
$$\Lambda\left(CH(\Omega)\right) \subseteq \Lambda\left(CH\bigl(cl(\Omega)\bigr)\right)=cl(\Omega)~,$$
where the last equality is achieved using Lemma \ref{convexhulls}. Finally, $\Lambda\left(CH(\Omega)\right)$ is
a closed subset of $\partial G$ containing $\Omega$ (by part (b) of Lemma \ref{limitsets}), which implies
statement of the Remark \ref{LCH-sets}.

The following lemma (in a slightly different form) can be found in \cite[Cor. to Lemma 13]{Swenson}:
\begin{lemma} \label{stabilizer} Suppose $\Omega \subset \partial G$ is a  subset having at least
two distinct points. Then $\Lambda\bigl(St_G(\Omega)\bigr) \subseteq cl(\Omega)$.
\end{lemma}

\underline{Proof.} Since $\Omega$ has at least two points, it makes sense to consider the convex hull $CH(\Omega)$.
Observe that for any $g \in St_G(\Omega)$, $gCH(\Omega) \subseteq CH(\Omega)$: the left
translation by the element $g \in G$ is an isometry of \ga, therefore a bi-infinite geodesic $(\alpha,\beta)$,
$\alpha,\beta \in \Omega$ goes to a bi-infinite geodesic $(g\circ\alpha,g\circ\beta) \subset CH(\Omega)$ since
$\Omega$ is $St_G(\Omega)$-invariant.

Fix any point $x \in CH(\Omega)$. By our observation above, $St_G(\Omega)x \subset CH(\Omega)$, hence
$\Lambda\bigl(St_G(\Omega)x\bigr) \subset \Lambda\bigl(CH(\Omega)\bigr)$. The claim of the lemma now follows
by applying Lemma \ref{limitsets}.(d) and Remark \ref{LCH-sets}. $\square$

\section{\large Known Results and Examples} \label{known}
In this section we list some known results about the boundary and limit sets and give (counter)examples
in order to motivate the rest of the paper where we extend these results to larger classes of subsets.


\vspace{.15cm}
{\bf Result 1:} if $A$ and $B$ are quasiconvex subgroups of a hyperbolic group $G$ then
$A \approx B$ if and only if $\Lambda(A)=\Lambda(B)$.

\vspace{.15cm}
Indeed, the necessity follows by Lemma \ref{limitsets}.(e). For proving the sufficiency we note that
by \cite[thm. 8]{Swenson}
$\Lambda (A \cap B) = \Lambda(A) \cap \Lambda(B)=\Lambda (A) = \Lambda (B)~.$ But
$A\cap B \le A$ and $A \cap B \le B$, so, by Lemma \ref{intersection} and \cite[thm. 4]{Swenson}
$$|A:(A \cap B)| < \infty,~|B:(A \cap B)| < \infty,$$ i.e. the subgroups $A$ and $B$ are commensurable.
Hence $A \approx B$.

However, if one removes at least one of the conditions on $A$ and $B$, the claim of the result 1 fails:

\vspace{.15cm}
{\bf \underline {Example 1.}} Let $G=F(x,y)$ -- the free group with two free generators $x,y$. Define
$A =\{x^n~|~n \ge 0\}$, $B=\{x^ny^m~|~0 \le m \le n\}$. These are quasiconvex subsets (not subgroups) of $G$
because any prefix of an element from one of these sets is still contained
in the same set. Evidently, $\Lambda (A) = \{x^\infty\}$. \\
Suppose $\left(x^{n_i}y^{m_i}\right)_{i \in \N}$, $0 \le m_i \le n_i$, $i \in \N$,
is a sequence converging to infinity in $B$. If the sequence of integers $(n_i)_{i \in \N}$ is bounded then
the sequence $(m_i)_{i \in \N}$ is also bounded, hence the set the group of elements in
$\left(x^{n_i}y^{m_i}\right)_{i \in \N}$ is finite which contradicts to the definition of a sequence that converges
to infinity. Thus, $\displaystyle \sup_{i\in \N} \{n_i\} = \infty$ and, passing to a subsequence, we can assume
$\displaystyle \lim_{i\to \infty} n_i = \infty$. Then
$$(x^{n_i}y^{m_i}|x^{n_i})_{1_G} = n_i \to \infty~\mbox{as}~i \to \infty~,$$
thus, $\displaystyle \lim_{i \to \infty} \left(x^{n_i}y^{m_i} \right) = \lim_{i \to \infty} x^{n_i} = x^\infty$.
Therefore, $\Lambda (B)=\{x^\infty\} = \Lambda(A)$  but $A \not \approx B$.

\vspace{.15cm}
{\bf \underline {Example 2.}} If $G$ is an arbitrary hyperbolic group and $H$ is its infinite normal subgroup
of infinite index, then $H$ is not quasiconvex by Corollary \ref{cor2} and $\Lambda(H) = \Lambda(G) = \partial G$
by \cite[Lemma 3.8.(2)]{K-S}, thus the quasiconvexity of $A,B$ in result 1 is important.

\vspace{.15cm}
{\bf Result 2:} If $A$ is a quasiconvex subgroup of a hyperbolic group $G$ then we have an equality
in the Remark \ref{rem4}: $Comm_G(A) = St_G \bigl(\Lambda(A) \bigr)$.

\vspace{.15cm}
By \cite[thm. 17]{Swenson} or \cite[cor. 3.10]{K-S}, $VN_G(A) = St_G \bigl(\Lambda(A) \bigr)$
and $VN_G(A) = Comm_G(A)$ by Remark \ref{rem1}.

Again, both of the conditions of $A$ being a subgroup and $A$ being quasiconvex are not redundant:

\vspace{.15cm}
{\bf \underline {Example 3.}} Choose $G=F(x,y)$ as in example 1 and let
$$B= \{x^ny^m~|~0 \le m \le n^2, n\ge 0\},~C=\{x^{-n}~|~n \in \N\},~A= B\cup C~.$$
$A$ is quasiconvex since $B$ and $C$ are, $\Lambda(A) = \Lambda(B) \cup \Lambda(C)=\{x^\infty,x^{-\infty}\}$
($\Lambda B =\{x^\infty\}$ by a similar argument to the one presented in example 1).
Then $St_G \bigl(\Lambda(A) \bigr) = \langle x \rangle$ -- the infinite cyclic subgroup generated
by $x$.\\
Let's show that $Comm_G(A) = \{1_G\}$. By Remark \ref{rem4} and since $Comm_G(A)$ is a subgroup, it is enough to prove that
$x^{-k} \notin Comm_G(A)$ for any integer $k>0$. Indeed, for any $n > k$ ~~$x^{n-k}y^{n^2} \in x^{-k}A$ and
$$d\left(x^{n-k}y^{n^2},A\right)=d\left(x^{n-k}y^{n^2},x^{n-k}y^{(n-k)^2}\right)=n^2-(n-k)^2=2nk-k^2 \to
\infty$$ when $n \to \infty$. Implying that $x^{-k}A \not \approx A$.

\vspace{.15cm}
{\bf \underline {Example 4.}} Consider a finitely generated group $M$ containing a normal subgroup
$N \lhd M$ and an infinite subnormal subgroup $K \lhd N$ such that $|M:N|=\infty$, $|N:K|=\infty$ and for any
$x \in M \backslash N$~~ $xKx^{-1} \cap K = \{1_M\}$ (for example, one can take $M = \Z~wr~\Z$).\\
Then $M$ is isomorphic to a quotient of some free group $G$ of finite rank by its normal subgroup
$H$: $M \cong G/H$. Let $\phi: G \to G/H$ be the natural homomorphism and $A,B \le G$ be the preimages
of $K$ and $N$ under $\phi$ correspondingly. Then $H \lhd A \lhd B \lhd G$, $|G:A|=\infty$.
$\Lambda (A) =\Lambda(B) = \Lambda(G)=\partial G$ by
\cite[lemma 3.8,(2)]{K-S}, hence $St_G\bigl( \Lambda (A) \bigr) = G$. \\
We claim that $Comm_G(A) = B$. As we know $Comm_G(A)=VN_G(A)$, therefore $B \subset Comm_G(A)$. Now, for an
arbitrary $g \in G \backslash B$, by construction, one has $\phi(A \cap gAg^{-1})=\{1_M\}$, hence
$(A \cap gAg^{-1}) \subset H$. Since $K$ is infinite, we get $|A:H|=\infty$, and thus, $|A:(A \cap A^g)| = \infty$,
so, $g \notin Comm_G(A)$.

In this example the subgroup $A$ of $G$ is not quasiconvex and \\$|St_G\bigl(\Lambda(A)\bigr):Comm_G(A)|=\infty$.

\vspace{.15cm}
{\bf Result 3:} (\cite[Thm. 2]{Arzh},\cite[Lemma 3.9]{K-S}) If $A$ is an infinite quasiconvex subgroup of a
hyperbolic group $G$ then $A$ has a finite index in its commensurator $Comm_G(A)$.

\vspace{.15cm}
By Lemma \ref{sgpcompar}, the condition $|Comm_G(A):A|<\infty$ is equivalent to $$Comm_G(A) \preceq A~.$$
It is easy to construct an example of a quasiconvex subset (not subgroup) $A$ with exactly one limit point
demonstrating that the latter fails, more precisely, $Comm_G(H)$ can have two limit points.

However, in the next section this result will be extended to the class of all quasiconvex subsets $A$ with
$card(\Lambda(A)) \ge 2$.

\vspace{.15cm}
{\bf Result 4:} Let $A$ be a quasiconvex subgroup of a hyperbolic group $G$. Then $Comm_G(A)$ is quasiconvex.

\vspace{.15cm}
If the subgroup $A$ is infinite, this is a consequence of the result 3 by Remark \ref{rem2}. On the other hand,
if $A$ is finite, then $Comm_G(A)=G$.

Below we give an example of an infinite quasiconvex set $A \subset G$ such that $Comm_G(A)$ is not quasiconvex.

\vspace{.15cm}
{\bf \underline {Example 5.}} We use the counterexample $12$ from \cite{Swenson}. Again, let $G=F(x,y)$ be the free
group of rank 2. Let $K = \langle x^nyx^{-n}~|~n \ge 0 \rangle$. It is shown in \cite{Swenson}
that $\Lambda(K)$ is not a limit set of a quasiconvex subgroup in $G$ (because $\Lambda (K)$ is not "symmetric":
$x^\infty \in \Lambda (K)$ but $x^{-\infty} \notin \Lambda(K)$). \\
As the subgroup $K$ is infinite, we can consider the convex hull $A=CH\bigl(\Lambda(K)\bigr)$.
By Lemma \ref{convexhulls} $A$ is quasiconvex and $\Lambda(A) = \Lambda(K)$ ( $\Lambda(K) \subset \partial G$
is closed by the claim (b) of Lemma \ref{limitsets} ). $K \subset Comm_G(K) \subset St_G\bigl(\Lambda(K)\bigr)$, hence,
as we saw in the proof of Lemma \ref{stabilizer}, $A$ is $K$-invariant. Consequently, $K \subset Comm_G(A)$.
Remark \ref{rem4} and Lemma \ref{stabilizer} imply
$$\Lambda(K) \subset \Lambda\bigl(Comm_G(A)\bigr) \subset \Lambda\left(St_G\bigl(\Lambda(A)\bigr)\right)
\subset \Lambda(A)=\Lambda(K)~.$$
Thus $\Lambda\bigl(Comm_G(A)\bigr) = \Lambda\left(St_G\bigl(\Lambda(A)\bigr)\right)=\Lambda(K)$, therefore the subgroups
\\$Comm_G(A)$ and $St_G\bigl(\Lambda(A)\bigr)$ are not quasiconvex.

\vspace{.15cm}
In the next section we are going to extend the results 1-3 to a broader class of quasiconvex subsets of the
hyperbolic group $G$. In particular, we will substitute the requirement for $A$ and $B$ to be subgroups by a
weaker condition.

\section{\large Tame Subsets} \label{generalization}
Again, let $G$ be a $\delta$-hyperbolic group with fixed finite generating set $\cal A$.

\vspace{.15cm}
{\bf \underline{Definition.}} A subset $A$ of the group $G$ will be called {\it tame} if $A$ has
at least two limit points on $\partial G$ and $A \preceq CH\bigl(\Lambda(A)\bigr)$. I.e. there exists
$\nu \ge 0$ such that $A \subset {\cal O}_\nu(C)$ where $C=CH\bigl(\Lambda(A)\bigr)$.

\vspace{.15cm}
In particular, this definition implies that any tame subset is infinite.

\vspace{.15cm}
\remark \label{rem5} If $A$ and $D$ are subsets of $G$ such that $A \approx D$ and $A$ is tame then $D$
is also tame.

\vspace{.15cm}
Indeed, By Lemma \ref{limitsets}.(e) $\Lambda(A)=\Lambda(D)$, hence,
$$D \preceq A \preceq CH\bigl(\Lambda(A)\bigr)=CH\bigl(\Lambda(D)\bigr)~.$$

Thus "tameness" of a subset is preserved under the equivalence relation "$\approx$".

\begin{lemma} \label{tamesets} Let $A,B,C,D$ be non-empty subsets of the group $G$ where $A$ and $B$ are tame,
$C$ is finite and $D$ is arbitrary. Let $H \le G$ be an infinite subgroup. Then the following sets are tame:
1) $A \cup B$; 2) $A\cup C$; 3)$A\cdot C$; 4) $D\cdot A$; 5) $H$.
\end{lemma}

\underline{Proof.} 1) Since $\Lambda(A),\Lambda(B) \subset \Lambda(A \cup B)$, we have
$$CH\bigl(\Lambda(A)\bigr) \cup CH\bigl(\Lambda(B)\bigr) \subseteq CH\bigl(\Lambda(A \cup B)\bigr)~.$$
$A \preceq CH\bigl(\Lambda(A)\bigr)$ and $B \preceq CH\bigl(\Lambda(B)\bigr)$ by conditions of the lemma, hence
$$A \cup B \preceq CH\bigl(\Lambda(A)\bigr) \cup CH\bigl(\Lambda(B)\bigr) \preceq CH\bigl(\Lambda(A \cup B)\bigr)~,$$
which shows that $A\cup B$ is tame.

2) and 3) are immediate consequences of the fact that $A \cup C \approx A$ , $A \cdot C \approx A$, and Remark \ref{rem5}.

4) Denote $K=CH\bigl(\Lambda(A)\bigr)$. By definition, $A \preceq K$, therefore, $DA \preceq DK$.
Now, since for every $y \in D$ $yK = CH\bigl(\Lambda(yA)\bigr) \subset CH\bigl(\Lambda(DA)\bigr)$,
we obtain $DK \subset CH\bigl(\Lambda(DA)\bigr)$. Hence $DA \preceq CH\bigl(\Lambda(DA)\bigr)$.

5) The set $CH\bigl(\Lambda(H)\bigr)$ is $H$-invariant by Lemma \ref{invariant}, therefore, for any
\\ $x \in CH\bigl(\Lambda(H)\bigr)$ we have $Hx \subset CH\bigl(\Lambda(H)\bigr)$. But $H \preceq Hx$,
hence $H$ is a tame subset. $\square$

\vspace{.15cm}
In particular, this lemma shows that any infinite set $U$ that is a finite union of quasiconvex products
in a hyperbolic group $G$ is tame.

In example 3 from section \ref{known} we constructed a quasiconvex subset $A$ in the group $G=F(x,y)$ with
exactly two limit points
$x^\infty,x^{-\infty}$. Therefore, $CH\bigl(\Lambda(A)\bigr)$ consists of one bi-infinite geodesic and
$CH\bigl(\Lambda(A)\bigr) \cap G = \{x^n~|~n \in \Z\}$. Now, for each $n\in \N$, $x^ny^{n^2} \in A$ and
$d\left(x^ny^{n^2},CH\bigl(\Lambda(A)\bigr)\right)=n^2 \to \infty$, as $n \to \infty$. Thus, the subset $A$
from example 3 is not tame.

\begin{lemma} \label{main-tame} Suppose $A$ is a tame subset of a hyperbolic group $G$ and $B \subseteq G$ is a
quasiconvex subset such that $\Lambda(A) \subseteq \Lambda(B)$. Then $A \preceq B$.
\end{lemma}

\underline{Proof.} By the conditions of the lemma, $A \preceq CH\bigl(\Lambda(A)\bigr) \preceq CH\bigl(\Lambda(B)\bigr)$.
Therefore, it remains to show that $CH\bigl(\Lambda(B)\bigr) \preceq B$, i.e. there exists $\varkappa \ge 0$
such that $CH\bigl(\Lambda(B)\bigr) \subset {\cal O}_\varkappa(B)$.

Let $\varepsilon$ be the quasiconvexity constant for $B$. Consider any $x \in CH\bigl(\Lambda(B)\bigr)$. By
definition, there exist $\alpha,\beta \in \Lambda(B)$ such that $x \in (\alpha,\beta)$.\\
Let $r_1,r_2: [0,\infty) \to \Gamma(G,{\cal A})$ be the geodesic half-lines obtained by bisecting $(\alpha,\beta)$
at the point $x$. Thus, $r_1(0)=r_2(0)=x$, $\displaystyle \lim_{i \to \infty} r_1(i)=\alpha$,
 $\displaystyle \lim_{i \to \infty} r_2(i)=\beta$ (see Figure 4).

\begin{figure}
\begin{center}
\input{pic4.tex}
\end{center}
\begin{center}
Figure 4
\end{center}

\end{figure}

There are sequences $(a_i)_{i \in \N}$ and $(b_i)_{i \in \N}$ in $B$ converging to infinity such that
$\displaystyle \lim_{i \to \infty} a_i=\alpha$, $\displaystyle \lim_{i \to \infty} b_i=\beta$.
Hence, $(r_1(i)|a_i)_x \to \infty$, $(r_2(i)|b_i)_x \to \infty$ as $i \to \infty$. Consequently,
for some $n \in \N$ we have
$$(r_1(n)|a_n)_x > 2\delta,~~(r_2(n)|b_n)_x > 2\delta~.\eqno (*)$$

\vspace{.15cm}
\remark \label{rem6} Let $PQR$ be a geodesic triangle in \ga and $(P|Q)_R>2\delta$. Then $d(R,[P,Q]) > 2\delta$.

\vspace{.15cm}
Indeed, assume, by contradiction, that there exists $S \in [P,Q]$ satisfying $d(R,S) \le 2\delta$. By definition of
the Gromov product,
$$(P|Q)_R=\frac12 \bigl(d(P,R)+d(Q,R)-d(P,Q) \bigr) \le $$
$$ \le \frac12 \bigl( d(P,S)+d(S,R)+d(Q,S)+d(S,R) -d(P,Q)\bigr)~.$$
But $d(P,S)+d(Q,S)=d(P,Q)$ since $[P,Q]$ is a geodesic segment, therefore
$(P|Q)_R \le d(S,R) \le 2\delta$. A contradiction.

\vspace{.15cm}
Consider, now, the geodesic quadrangle in \ga with vertices $a_n$, $r_1(n)$, $r_2(n)$, $b_n$.
$x \in [r_1(n),r_2(n)]$ (Fig. 4). Applying ($*$) and the Remark \ref{rem6} we obtain
$$d\bigl(x,[a_n,r_1(n)]\bigr) > 2\delta,~d\bigl(x,[b_n,r_2(n)]\bigr)>2\delta~.$$

Since the Cayley graph \ga is $\delta$-hyperbolic, all quadrangles are $2\delta$-slim, thus
$$[r_1(n),r_2(n)] \subset {\cal O}_{2\delta}\bigl([a_n,r_1(n)] \cup [b_n,r_2(n)] \cup [a_n,b_n] \bigr)~.$$
Consequently, $d(x,[a_n,b_n]) \le 2\delta$. $a_n,b_n \in B$ and $B$ is $\varepsilon$-quasiconvex, therefore
$[a_n,b_n] \subset {\cal O}_\varepsilon(B)$.

So, $d(x,B) \le 2\delta + \varepsilon$ for every $x \in CH\bigl(\Lambda(B)\bigr)$. After denoting
$\varkappa = 2\delta+\varepsilon$ we
achieve $CH\bigl(\Lambda(B)\bigr) \subset {\cal O}_\varkappa (B)$. The Lemma \ref{main-tame} is proved. $\square$

\begin{cor} \label{cor4}
Let $A \preceq B$ be subsets of $G$ where $A$ has at least two limit points on $\partial G$ and
$B$ is quasiconvex. Then $Comm_G(A) \preceq B$, $St_G \bigl( \Lambda(A)\bigr) \preceq B$.
\end{cor}

\underline{Proof.} By Remark \ref{rem4} it is enough to prove the second inequality. Using lemmas \ref{stabilizer},
\ref{limitsets} we get
$$\Lambda \Bigl(St_G\bigl(\Lambda(A)\bigr)\Bigr) \subseteq \Lambda(A) \subseteq \Lambda(B)~.$$
Since any subgroup is a tame subset (Lemma \ref{tamesets}), applying Lemma \ref{main-tame} we obtain
$St_G\bigl(\Lambda(A)\bigr) \preceq B$. $\square$

\vspace{.15cm}
The statement of corollary \ref{cor2} can be generalized as follows:

\begin{cor} \label{cor5} Let $G$ be a hyperbolic group and let $U$ be a finite union of quasiconvex products of
infinite index in $G$. Suppose $A \subset U$ is an infinite subset. Then $|G:Comm_G(A)|=\infty$.
\end{cor}

\underline{Proof.} First, notice that since $A$ is infinite, $U$ is also infinite, hence
$U$ has an infinite member $H \le G$ of infinite index. Therefore, $G$ is non-elementary.
There are two possibilities: either $card\bigl(\Lambda(A)\bigr) \le 1$ or $card\bigl(\Lambda(A)\bigr) \ge 2$.

In the first case $card\bigl(\Lambda(A)\bigr) = 1$ (by Lemma \ref{limitsets}.(a)), so
$\Lambda(A) = \{\alpha\} \in \partial G$. Hence $St_G\bigl(\Lambda(A)\bigr)=St_G(\{\alpha\})$ is an elementary
subgroup, thus $Comm_G(A)$ is also elementary. Consequently, $|G:Comm_G(A)|=\infty$.

In the second case, when $card\bigl(\Lambda(A)\bigr) \ge 2$, we can use corollary \ref{cor4} to obtain
$Comm_G(A) \preceq U$. But if $Comm_G(A)$ were quasidense (or, equivalently, \\ $|G:Comm_G(A)|<\infty$),
 we would have $G \preceq Comm_G(A) \preceq U$, so
$U$ would have to be also quasidense. The latter contradicts to the statement of corollary \ref{non-quasidense}. $\square$

\vspace{.15cm}
Now we are going to extend the results 1,2 to all tame quasiconvex subsets:
\begin{prop} \label{b-rigid} Suppose $A$ and $B$ are tame quasiconvex subsets of a hyperbolic group $G$.
Then $A \approx B$ if and only if $\Lambda(A)=\Lambda(B)$.
\end{prop}

\underline{Proof.} The necessity is given by Lemma \ref{limitsets}.(e); the sufficiency immediately
follows from Lemma \ref{main-tame}. $\square$

\begin{cor}
Let $A$ be a subset of a hyperbolic group $G$
having at least two distinct limit points on $\partial G$. Then
the following two conditions are equivalent:

1) $A$ is tame and quasiconvex;

2) $A \approx CH\bigl(\Lambda(A)\bigr)$.
\end{cor}

\underline{Proof.}  Denote $C=CH\bigl(\Lambda(A)\bigr)$. By Lemma \ref{limitsets}.(b) $
\Lambda(A)$ is a closed subset of $\partial G$, hence from Lemma \ref{convexhulls}.(b) we get
$\Lambda(A) = \Lambda(C)$.  Therefore $CH\bigl(\Lambda(C)\bigr)=C$ implying that $C$ is tame.
$C$ is quasiconvex by Lemma \ref{convexhulls}.(a).

Now, 2) follows from 1) by Proposition \ref{b-rigid}.

1) follows from 2) because the equivalence relation "$\approx$" preserves quasiconvexity and "tameness"
of a subset. $\square$

\begin{prop} \label{commeqstab} For any tame quasiconvex subset $A$ of a hyperbolic group $G$
$Comm_G(A) = St_G\bigl(\Lambda(A)\bigr)$.
\end{prop}

\underline{Proof.} By Remark \ref{rem4} it is enough to show that $ St_G\bigl(\Lambda(A)\bigr) \subseteq Comm_G(A)$.
Take an arbitrary $g \in St_G\bigl(\Lambda(A)\bigr)$. Then $\Lambda(gA)=g\circ\Lambda(A)=\Lambda(A)$.
The subset $gA$ is tame and quasiconvex since $A$ is so, hence, by Proposition \ref{b-rigid}, $gA \approx A$. Thus,
$g \in Comm_G(A)$. Q.e.d. $\square$

\begin{prop} \label{result3} Let $G$ be a hyperbolic group and let $A \subset G$ be a quasiconvex subset
that has at least two distinct limit points on the boundary $\partial G$. Then
$St_G\bigl(\Lambda(A)\bigr) \preceq A$. Consequently, $Comm_G(A) \preceq A$.
\end{prop}

\underline{Proof.} Denote $H=St_G\bigl(\Lambda(A)\bigr)$. Using lemmas \ref{limitsets}.(b) and \ref{stabilizer}
we obtain $\Lambda(H) \subseteq \Lambda(A)$. $H\le G$ is a subgroup, hence it is tame (Lemma \ref{tamesets});
$A$ is quasiconvex by the conditions. The statement of the proposition now follows from Lemma \ref{main-tame}.
 $\square$

The Proposition \ref{result3} generalizes the result 3 from section \ref{known}.

Now, it is easy to see that the set $A$ from example 5 in section \ref{known} is tame and quasiconvex,
thus, $Comm_G(A) = St_G\bigl(\Lambda(A)\bigr)$. However $Comm_G(A)$ is not quasiconvex.
Thus we can not extend the result 4 from section \ref{known} in the same way we did the previous ones.

\section{Proof of the Theorem \ref{conjugacythm}} \label{conjsection}
A subset $B$ of the group $G$ will be called {\it normal} if for every $g \in G$~~ $B^g=B$.

\begin{lemma} \label{quasidense} Suppose $G$ is a non-elementary hyperbolic group and $A$ is its subset containing
an infinite normal subset $B$. Then $A$ is quasiconvex if and only if it is quasidense in $G$.
\end{lemma}

\underline{Proof.} The sufficiency is trivial (see part 3) of Remark \ref{rem2.5}). To show the necessity,
first we observe that $Comm_G(B) = G$, hence $St_G\bigl(\Lambda(B)\bigr)=G$ by Remark \ref{rem4}. Since $B$ is infinite $\Lambda(B)$ is a non-empty subset of the boundary $\partial G$, thus there exists
$\alpha \in \Lambda(B)$. Now, if $\Lambda(B)=\{\alpha\}$ then, as we know, $G=St_G\bigl(\Lambda(B)\bigr)=St_G(\{\alpha\})$
is an elementary subgroup of $G$, but $G$ is non-elementary.
Therefore $\Lambda(B)$ has at least two distinct points. Now we can apply Lemma \ref{stabilizer} to obtain
$$\partial G =\Lambda(G) = \Lambda\left( St_G\bigl(\Lambda(B)\bigr) \right) \subseteq \Lambda(B) \subseteq \partial G$$
and conclude that $\Lambda(B) = \partial G$. $B \subset A$, thus $\Lambda(A)=\partial G$.
Finally, since $G$ is a tame subset of itself (by Lemma \ref{tamesets}) and $A$ is quasiconvex, we apply Lemma \ref{main-tame} and achieve $G \preceq A$. $\square$

\begin{prop} \label{prop5} Let $G$ be a hyperbolic group, $K,H_1,\dots,H_s$ be its subgroups,
where $K$ is non-elementary and $H_j$ are $\varepsilon_j$-quasiconvex for all $j=1,2,\dots,s$.
If $\displaystyle K\preceq \bigcup_{j=1}^s H_j^G$ then
for some $k \in \{1,\dots,s\}$ and $g \in G$ one has $K \preceq H_k^g$, i.e. the index $|K:(K\cap H_k^g)|$ is finite.
\end{prop}

\underline{Proof.} Denote $\varepsilon = max\{\varepsilon_1,\dots,\varepsilon_s\}$ and assume the contrary to
the statement we need to prove.
Then by Theorem\ref{prop1} there exists an element $x \in K$ of
infinite order satisfying $\displaystyle \langle x \rangle \cap \bigcup_{j=1}^s H_j^G = \{1_G\}$.
Now, since $K$ is non-elementary we can apply Proposition \ref{thm1} to obtain an element $y \in K$
with $o(y)=\infty$ such that $\langle x \rangle^G \cap \langle y \rangle $ is trivial. Consequently,
$\langle x \rangle \cap \langle y \rangle^G=\{1_G\}$.

By the conditions of the proposition
$$ K \subset  \bigcup_{i=1}^r \left( \bigcup_{j=1}^s H_j^G\right)b_i =
\bigcup_{i=1}^r \left( \bigcup_{j=1}^s H_j^Gb_i\right)$$
for some $b_1,b_2,\dots,b_r \in G$. For every $n \in \N$~~ $x^ny^n \in K$, hence there
exist $j_0 \in \{1,\dots,s\}$ and $i_0 \in \{1,\dots,r\}$ such that $x^ny^n \in (H_{j_0})^G \cdot b_{i_0}$ for
infinitely many $n \in \N$. Without loss of generality, assume $j_0=1$, $i_0=1$.

Let $w_1,w_2,w_3$ be words in alphabet $\cal A$ representing elements $x,y,b_1^{-1}$ correspondingly.


Observe that the word $(w_1^nw_2^nw_3)^l$ represents the element $(x^ny^nb_1^{-1})^l$. For infinitely many
$n \in \N$, $n \ge N$, there exists $g_n \in G$ such that $x^ny^nb_1^{-1} \in g_nH_1g_n^{-1}$, thus,
$(x^ny^nb_1^{-1})^l \in g_nH_1g_n^{-1}$ for all $l \in \N$.

For any such $n$  and $l$, define the points $Y_j$
in \ga, $j=0,1,\dots,2l$, as follows: $Y_0=g_n$, $Y_1=g_nx^n$, $Y_2=g_nx^ny^nb_1^{-1}$, $Y_3=g_nx^ny^nb_1^{-1}x^n$,
$Y_4=g_n(x^ny^nb_1^{-1})^2$, \dots, $Y_{2l-1}=g_n(x^ny^nb_1^{-1})^{l-1}x^n$, $Y_{2l}=g_n(x^ny^nb_1^{-1})^l$.
Consider the path $q$ starting at $Y_0$ and labelled by the word $(w_1^nw_2^nw_3)^l$. Thus, $q$ ends at
$Y_{2l}$ and passes through each $Y_j$, $j=1,\dots,2l-1$. By lemmas \ref{power} and \ref{close}
applied to the segment of $q$ between $Y_j$ and $Y_{j+1}$ and the geodesic $[Y_j,Y_{j+1}]$, for each $j=0,\dots,2l-1$,
there exists a constant $\nu = \nu(w_1,w_2,w_3) \ge 0$ such that $q \subset {\cal O}_\nu([Y_0,\dots, Y_{2l}])$.

By the choice of $y$ and Lemma \ref{quasintersection} there exists $\varkappa=\varkappa(x,y,b_1) \ge 0$ such
that $(x^{-n}|y^n)_{1_G} \le \varkappa$ ~and~ $(b_1y^{-n}b_1^{-1}|x^n)_{1_G} \le \varkappa$ ~for all $n$.

Denote $C_0=\varkappa + |b_1|_G$. Then, for any odd $j \in \{1,2,\dots,2l-1\}$,
$$(Y_{j-1}|Y_{j+1})_{Y_j} = (x^{-n}|y^nb_1^{-1})_{1_G} \le (x^{-n}|y^n)_{1_G} + |b_1|_G \le C_0.$$
And for any even $j \in \{1,2,\dots,2l-1\}$,
$$(Y_{j-1}|Y_{j+1})_{Y_j} = (b_1y^{-n}|x^n)_{1_G} \le (b_1y^{-n}b_1^{-1}|x^n)_{1_G} + |b_1|_G \le C_0.$$

Choose an arbitrary $C_1 > 12(C_0+\delta)$. Lemma \ref{power} implies that for any sufficiently large $n$,
$\|[Y_{j},Y_{j+1}]\|>C_1$, $j=1,2,\dots,2l-1$, hence, we can use lemmas \ref{brokenlines1} and \ref{brokenlines2}
to obtain
$$[Y_0,\dots,Y_{2l}] \subset {\cal O}_{2C_0}([Y_0,Y_{2l}])~\mbox{ and }~\|[Y_0,Y_{2l}]\| \ge
\|[Y_0,\dots,Y_{2l}]\|/2 \ge \frac{2lC_0}2~.$$

Let $k$ be an even number from the set $\{l-1,l\}$. Then the subpath $p$ of $q$ between $Y_k$
and $Y_{k+1}$ is labelled by the word $w_1^n$ and $p \subset {\cal O}_\nu([Y_{k}, Y_{k+1}])$. As we showed, there are
points $U,V \in [Y_0,Y_{2l}]$ such that $d(Y_{k},U) \le 2C_0$ and $d(Y_{k+1},V) \le 2C_0$. Since the quadrangles
in a $\delta$-hyperbolic space are $2\delta$-slim, we obtain $[Y_{k}, Y_{k+1}] \subset {\cal O}_{2C_0+2\delta}([U,V])$.
Applying lemmas lemmas \ref{brokenlines1} and \ref{brokenlines2} to $[Y_0,\dots,Y_k]$ and $[Y_{k+1},\dots,Y_{2l}]$
we achieve
$$d(Y_0,U) \ge d(Y_0,Y_k)-2C_0\ge \frac{\|[Y_0,\dots,Y_k]\|}{2} -2C_0 \ge \frac{C_1(l-1)}{2}-2C_0 ~\mbox{ and }$$
$$d(V,Y_{2l}) \ge d(Y_{k+1},Y_{2l})-2C_0\ge \frac{\|[Y_{k+1},\dots,Y_{2l}]\|}{2} -2C_0 \ge \frac{C_1l}{2}-2C_0~.$$

Hence, if $l\in \N$ is sufficiently large, we will have $d(Y_0,Y_{2l})>2|g_n|_G$, \\
$d(Y_0,U)>|g_n|_G$ and $d(V,Y_{2l})>|g_n|_G$ (and, similarly, $d(U,Y_{2l})>|g_n|_G$, $d(V,Y_0)>|g_n|_G$).
%

Let $X_1=1_G$, $X_2 = g_n(x^ny^nb_1^{-1})^lg_n^{-1} \in H_1$ for infinitely many $n \in \N$. By Lemma \ref{quadrangle}
applied to the geodesic quadrangle $X_0Y_0Y_{2l}X_1$, the subsegment $[U,V]$ of $[Y_0,Y_{2l}]$ lies in
the $2\delta$-neighborhood
of $[X_0,X_1]$ and $[X_0,X_1]\subset {\cal O}_{\varepsilon}(H_1)$ (for infinitely many $n$ and sufficiently large $l$).
Accumulating all of the above, we obtain
$$p \subset {\cal O}_{\nu+2C_0+2\delta+\varepsilon}(H_1)$$
(for infinitely many $n \in \N$ and sufficiently large $l =l(n)\in \N$). As in the proof of Proposition \ref{thm1},
the latter implies that for some $t \in \N$ we have $x^t \in H_1^G$ -- a contradiction to the construction of
$x$.

The proposition is proved. $\square$

\vspace{.15cm}
\underline{Proof of Theorem\ref{conjugacythm}.} It is given that $A= a_1^G \cup a_2^G \dots \cup a_s^G$
for some elements $a_1,\dots,a_s \in G$.

Suppose that $A$ is quasiconvex. $A$ is infinite, therefore by Lemma \ref{quasidense}
the subset $A$ is quasidense in $G$, hence
$$G \preceq A \subset \bigcup_{j=1}^s \langle a_j \rangle^G~.$$
But since $G$ is non-elementary, $|G:\langle a_j \rangle^g|=\infty$ for each $j=1,2,\dots,s$, and for any $g\in G$.
Thus, we achieve a contradiction with Proposition \ref{prop5} applied to $K=G$. Q.e.d. $\square$

\end{document}